\input amstex 
\documentstyle{amsppt} 
\magnification=1200 
 
\def\la{\lambda} 
\def\be{\beta} 
\def\al{\alpha} 
\def\La{\Lambda} 
\def\ka{\varkappa} 
\def\Lt{\Lambda^\theta} 
\def\om{\omega} 
\def\Om{\Omega} 
\def\Ga{\Gamma} 
\def\R{\Bbb R} 
\def\C{\Bbb C} 
\def\Z{\Bbb Z} 
\def\dimth{\varkappa_\theta} 
\define\pf{\operatorname{Pf}} 
 
\def\th{\theta} 
\def\Y{\Bbb Y} 
\def\J{\Bbb J(\theta)} 
\def\JJ{\Bbb J}
\def\K{\Bbb K} 
\def\G{\Bbb G} 
\def\SS{\Bbb S} 
\def\dimG{\dim_\G} 
\def\Harm{{\Cal H}} 
\def\P{\Cal P} 
\def\Pf{\operatorname{Pf}} 
 
\def\X{\frak X} 
\def\M{\Cal M}

\def\dhr{\downharpoonright} 
 
\def\tht{\thetag} 
\def\wt{\widetilde} 
 
\NoRunningHeads 
\TagsOnRight 
 
\topmatter 
\thanks
One of the authors (G.O.) was supported by the Russian Foundation for Basic Research under grant 98--01--00303.
\endthanks
\title Harmonic functions on multiplicative graphs and interpolation 
polynomials  
\endtitle 
 
\author Alexei Borodin and Grigori Olshanski 
\endauthor 
 
\abstract We construct examples of nonnegative harmonic functions on 
certain graded graphs: the Young lattice and its generalizations. Such functions first emerged in harmonic analysis on 
the infinite symmetric group. Our method relies on multivariate 
interpolation polynomials associated with Schur's S and P functions 
and with Jack symmetric functions. As a by--product, we compute 
certain Selberg--type integrals.  
\endabstract 
 
\toc 
\widestnumber\head{\S7.} 
\head \S0. Introduction \endhead 
\head \S1. The general formalism\endhead 
\head \S2. The Young graph \endhead 
\head \S3. The Jack graph \endhead  
\head \S4. The Kingman graph \endhead 
\head \S5. The Schur graph \endhead 
\head \S6. Finite--dimensional specializations \endhead 
\subhead 6.1. Truncated Young branching \endsubhead 
\subhead 6.2. $\Gamma$--shaped Young branching \endsubhead 
\subhead 6.3. Truncated Kingman branching \endsubhead 
\subhead 6.4. Truncated Schur branching \endsubhead 
\head \S7. Appendix \endhead 
\head{} References \endhead 
\endtoc 
 
\endtopmatter 
\document 
 
\head \S0. Introduction \endhead 
 
Let $\Y$ denote the lattice of Young diagrams ordered by inclusion. 
For $\mu,\la\in\Y$, we write $\la\searrow\mu$ if $\la$ covers $\mu$, 
i.e., $\la$ differs from $\mu$ by adding a box. We consider $\Y$ as a 
graph whose vertices are arbitrary Young diagrams $\mu$ and the edges 
are couples $(\mu,\la)$ such that $\la\searrow\mu$. We shall call $\Y$ the {\it Young graph}. 
A function $\varphi(\mu)$ is called a {\it harmonic function} on the Young 
graph \cite{VK} if it satisfies the condition  
$$ 
\varphi(\mu)=\sum_{\la:\, \la\searrow \mu} \varphi(\la), \qquad  
\forall\mu\in\Y. \tag0.1 
$$ 
We are interested in nonnegative harmonic functions $\varphi$ 
normalized at the empty diagram: $\varphi(\varnothing)=1$. Such  
functions form a convex set denoted as $\Harm^+_1(\Y)$. 
 
The functions $\varphi\in\Harm^+_1(\Y)$ have an important 
representation--theoretic meaning: they are in a natural bijective 
correspondence with central, positive definite, normalized functions 
on the infinite symmetric group $S(\infty)$, see \cite{VK}, 
\cite{KV2}. Thoma's description of characters on $S(\infty)$ means 
that the extreme points of $\Harm^+_1(\Y)$ form an 
infinite--dimensional simplex $\Om$ (called the Thoma simplex), see 
\cite{T}, \cite{VK}, \cite{KV2}, \cite{W}. For general elements 
$\varphi\in\Harm^+_1(\Y)$, there is a (unique) Poisson--type 
integral representation, 
$$ 
\varphi(\la)=\int_\Om K(\la,\om)\, P(d\om), \qquad \forall \la\in\Y, 
\tag0.2  
$$ 
where $P$ is a probability measure on $\Om$ (the `boundary measure' 
for $\varphi$) and $K(\la,\om)$ is a positive function on 
$\Y\times\Om$ (the `Poisson kernel' or `Martin kernel' for $\Y$), see 
\cite{KOO}. Note that any probability measure $P$ on $\Om$ gives rise 
to an element of $\Harm^+_1(\Y)$; in particular, the extreme 
$\varphi$'s are exactly the functions $K(\cdot,\om)$ corresponding to 
Dirac measures on $\Om$.   
 
This abstract result shows how large $\Harm^+_1(\Y)$ is but it does 
not explain how to construct explicitly nonextreme 
functions $\varphi$ or what nonextreme $\varphi$'s could be 
interesting for applications. 
 
Concrete examples of nonextreme functions $\varphi$ first emerged in 
\cite{KOV} in connection with a problem of harmonic analysis on the
infinite symmetric
group $S(\infty)$. These functions, denoted as $\varphi_{zz'}$, 
depend on two parameters, and the corresponding `boundary measures' 
$P_{zz'}$ govern the spectral decomposition of certain natural unitary 
representations. \footnote{The measures $P_{zz'}$ are very 
interesting objects. They are studied in detail in our papers 
\cite{P.I} -- \cite{P.V}, \cite{BO1}, \cite{BO2}.}  
 
The explicit expression of the functions $\varphi_{zz'}$ (see formula 
\tht{2.4} below) has an interesting combinatorial structure which 
raises a number of questions. For instance, one can ask whether 
there exist similar families of harmonic functions for other graphs. 
The answer is affirmative: \cite{B1}, \cite{Ke5}. \footnote{Another 
question, a characterization of the functions of type 
$\varphi_{zz'}$, was examined in \cite{B1}, \cite{Ro}.} 
 
The paper \cite{B1} concerns the graph $\SS$ of shifted Young 
diagrams which is related to projective representations of the 
symmetric groups.  
 
The paper \cite{Ke5} contains a generalization in another direction: 
a deformation of the family $\{\varphi_{zz'}\}$, which is consistent 
with a deformation of the basic equation \tht{0.1}: 
$$ 
\varphi(\mu)=\sum_{\la:\, \la\searrow\mu}\dimth(\mu,\la)\varphi(\la), 
\qquad \forall \mu. \tag0.3 
$$  
Here $\th>0$ is the deformation parameter and $\dimth(\mu,\la)>0$ are the 
coefficients that arise in (the simplest case of) Pieri's rule 
for Jack symmetric functions with parameter $\th$. The initial 
situation corresponds to the particular value $\th=1$, when Jack 
symmetric functions coincide with Schur's $S$--functions. 
 
Note that in the limit as $\th\to0$, the harmonicity condition 
\tht{0.3} essentially coincides with the relation which defines 
partition structures in the sense of Kingman \cite{Ki1}, \cite{Ki2}, 
while the two-parameter family of harmonic functions constructed in 
\cite{Ke5} degenerates to the famous Ewens partition structures 
\cite{Ew} and its generalization due to Pitman, see 
\cite{Pi}, \cite{PY}, \cite{Ke4}.   
 
In the present paper, we propose a simple combinatorial construction, 
which allows us to get, in a unified way, all these concrete examples of 
harmonic functions as well as some new ones. In the new examples, the 
`boundary measures' $P$ are supported by finite--dimensional simplices, and the Poisson integral representation leads to certain 
Selberg--type integrals. \footnote{A connection between Poisson 
integral representation of type \tht{0.2} and Selberg integrals was 
first exploited in \cite{Ke3}.} 
 
Our construction relies on the so--called shifted (or factorial) 
versions of Schur's $S$ and $P$ functions and of Jack symmetric  
functions. These new combinatorial functions arise in different 
topics, see, e.g., \cite{S}, \cite{KS}, \cite{OO}, \cite{OO2}, \cite{Ok1}, \cite{Ok2}. They are also called interpolation 
polynomials, because they give  solutions to certain  multivariate 
interpolation problems.  
 
The paper is organized as follows. In \S1, we expose the general 
formalism. In \S2, it is applied to the Young graph to derive the 
family $\{\varphi_{zz'}\}$. In \S\S3--5, we apply it to the Young 
graph with Jack edge multiplicities $\J$, next to the Kingman graph 
$\K$, and then to the Schur graph $\SS$; the arguments are quite 
similar. Section 6 is devoted to constructing harmonic functions 
of a different sort --- those with finite--dimensional `boundary 
measures'; here we also evaluate Selberg--type integrals. The 
final \S7 is an appendix on the Poisson integral representation.

\head \S1. The general formalism \endhead 
 
In this section, we deal with an abstract graph $\G$ satisfying 
certain conditions listed below. In the next sections, concrete 
examples of $\G$ will be considered. 
 
Our assumptions and conventions concerning $\G$ are as follows: 
 
$\bullet$ To simplify the notation, we identify the graph with its 
set of vertices. 
 
$\bullet$ The vertices are partitioned into levels, 
$\G=\G_0\sqcup\G_1\sqcup\G_2\sqcup\dots$, so that the endpoints of 
any edge lie on consecutive levels. That is, $\G$ is a graded graph.  
 
$\bullet$ The level of a vertex $\mu$ is denoted as $|\mu|$. If two 
vertices $\mu,\la$ form an edge, $|\la|=|\mu|+1$, then we write 
$\la\searrow\mu$ or $\mu\nearrow\la$.  
 
$\bullet$ All the levels $\G_n$ are finite. 
 
$\bullet$ The lowest level $\G_0$ consists of a single 
vertex denoted as $\varnothing$.  
 
$\bullet$ For any vertex $\mu$ there exists at least one vertex 
$\la\searrow\mu$ and for any vertex $\la\ne\varnothing$ there exists 
at least one vertex $\mu\nearrow\la$. This implies that the graph is 
connected.  
 
(Our main example is the Young graph, see \S2. ) 
 
$\bullet$ Finally, assume that we are given an {\it edge multiplicity 
function\/} which assigns to any edge $\mu\nearrow\la$ a strictly 
positive number $\ka(\mu,\la)$ --- its formal multiplicity. It should 
be emphasized that these numbers are not necessarily integers.  
 
(For the Young graph, all the formal multiplicities are equal to 1; 
graphs with nontrivial multiplicities are considered in \S3 and \S4.) 
 
A complex function $\varphi(\mu)$ on $\G$ is called a {\it harmonic function 
on the graph\/} $\G$ if it satisfies the relation 
$$ 
\varphi(\mu)=\sum_{\la:\, \la\searrow\mu} \ka(\mu,\la)\varphi(\la) 
\tag 1.1  
$$ 
for any vertex $\mu$ (the sum in the right--hand side is finite, because all the 
levels are finite). Let $\Harm(\G)$ denote the space of all harmonic 
functions endowed with the topology of pointwise convergence. Let 
$\Harm^+(\G)$ be the subset of nonnegative harmonic functions and 
$\Harm^+_1(\G)$ be the subset of the functions 
$\varphi\in\Harm^+(\G)$ with the normalization 
$\varphi(\varnothing)=1$.  
 
Clearly, $\Harm^+_1(\G)$ is a convex subset of $\Harm(\G)$. Moreover, 
it is a compact measurable space (here we again employ the 
finiteness assumption). We shall use some well--known general 
theorems about convex compact measurable sets which can be found, 
e.g., in \cite{Ph}. 
 
Let $\Om(\G)$ denote the set of extreme points in $\Harm^+_1(\G)$. 
This is a set of type $G_\delta$, hence, a Borel measurable set. Given 
$\om\in\Om(\G)$, let us denote by $K(\,\cdot\,, \om)$ the 
corresponding extreme harmonic function on $\G$. Note that  
$K(\mu,\,\cdot\,)$ is a Borel measurable function on $\Om(\G)$ for 
any fixed $\mu\in\G$.  
 
\proclaim{Theorem 1.1} For each element $\varphi\in\Harm^+_1(\G)$ 
there exists a unique probability measure $P$ on $\Omega(\G)$ such 
that  
$$ 
\varphi(\mu)=\int_{\Om(\G)} K(\mu,\om) P(d\om), \qquad 
\forall\mu\in\G. \tag 1.2 
$$ 
\endproclaim 
 
\demo{Proof} See \S7. \qed 
\enddemo 
 
We call \tht{1.2} the {\it Poisson integral representation\/} of the 
function $\varphi$. 
 
To any {\it path} $\tau$ going from a vertex $\mu$ to a vertex $\la$ with 
$|\la|>|\mu|$,  
$$ 
\tau=(\mu=\la_0\nearrow\la_1\nearrow\dots\nearrow\la_k=\la),  
\qquad k=|\la|-|\mu|, 
$$  
we assign its weight  
$$ 
w(\tau)=\prod_{i=1}^k\ka(\la_{i-1},\la_i) 
$$ 
and then set 
$$ 
\dimG(\mu,\la)=\sum_\tau w(\tau),  
\tag 1.3 
$$ 
summed over all paths from $\mu$ to $\la$. We extend this definition 
to all couples $(\mu,\la)$ by agreeing that $\dimG(\mu,\mu)=1$ and 
$\dim(\mu,\la)=0$ if $\mu\ne\la$ are such that there is no path from 
$\mu$ to $\la$. Next, we set $\dimG\la=\dimG(\varnothing,\la)$.  
 
In all examples of the graphs $\G$ considered in the present paper 
one can embed (the vertices of) $\G$ into $\Om(\G)$ in such a way that 
any point $\om\in\Om(\G)$ can be approximated by a sequence of vertices 
$\{\la(n)\in\G_n\}_{n=1,2,\dots}$, and for any such sequence 
$$ 
K(\mu,\om)=\lim_{n\to\infty}\frac{\dimG(\mu,\la(n))}{\dimG\la(n)}\,. 
$$ 
 
Given a function $\varphi\in\Harm^+_1(\G)$, we set for each $n$ 
$$ 
M_n(\la)=\dimG\la\cdot\varphi(\la), \qquad \la\in\G_n.  
\tag 1.4 
$$ 
Using the harmonicity relation \tht{1.1} and induction on $n$ one 
readily verifies that $\sum_{\la\in\G_n} M_n(\la)=1$. Thus, each $M_n$ is a 
probability distribution on $\G_n$.  
 
For all examples of the graphs $\G$ considered in this paper one can 
transfer the measure $M_n$ to $\Om(\G)$ via the 
embedding $\G\hookrightarrow\Om(\G)$ mentioned above. Then the 
measure $P$ appearing in the integral representation \tht{1.2} is the 
weak limit of the measures $M_n$ as $n\to\infty$.   
 
We say that $(\G,\ka(\,\cdot\,,\,\cdot\,))$ is a {\it multiplicative 
graph\/} \cite{KV1}, \cite{KV2}, if the following conditions are satisfied. First, 
the 1st floor $\G_1$ consists of a single vertex denoted by the 
symbol ``$(1)$''.  
Next, there exists a graded commutative unital algebra $A$ over $\R$, 
$A=A_0+A_1+\dots$, and a homogeneous basis $\{\P_\mu\}$ in $A$ indexed 
by the vertices $\mu\in\G$, such that $\P_\varnothing=1$ and $\deg 
\P_\mu=|\mu|$. Finally, for any $\mu$, 
$$ 
\P_\mu \P_{(1)}=\sum_{\la:\,\la\searrow\mu} \ka(\mu,\la)\P_\la\,.  
\tag 1.5 
$$ 
Note that this implies $\ka(\varnothing, (1))=1$. (All the graphs 
considered in the present paper are multiplicative. For instance, in 
the case of the Young graph, the algebra $A$ is the algebra of 
symmetric functions and the basis $\{\P_\mu\}$ is formed by the Schur 
functions.) 
 
Iterating the 
relation \tht{1.5} we get the expansion 
$$ 
\P_{(1)}^n=\sum_{\la\in\G_n}\dimG\la\cdot\P_\la\,, 
\tag 1.6 
$$ 
which is a useful tool for computing the dimensions $\dimG\la$. More 
generally, given $\mu\in\G_m$ and $n>m$, 
$$ 
\P_\mu\P_{(1)}^{n-m}=\sum_{\la\in\G_n}\dimG(\mu,\la)\cdot\P_\la\,. 
\tag 1.7 
$$ 
 
\proclaim{Theorem 1.2 \cite{KV1}} Let $\G$ be a 
multiplicative graph and let $A$ be the  
corresponding algebra. Given $\varphi\in\Harm^+_1(\G)$, let 
$\pi:A\to\C$ be the linear functional sending each $\P_\mu$ to 
$\varphi(\mu)$.  
Then $\varphi$ is extreme if and only if $\pi$ is multiplicative. 
\endproclaim 
 
Note that a linear functional $\pi:A\to\C$ corresponds to a function 
$\varphi\in\Harm^+_1(\G)$ if and only if $\pi(1)=1$, $\pi(\P_\mu)\ge0$ 
for any $\mu$, and $\pi$ factors through $A/(\P_{(1)}-1)A$.  
 
Now we shall explain our method of producing harmonic functions. 
Assume $A^*$ is a commutative algebra\footnote{The superscript $\ast$ does not 
mean the passage to a dual space.}, $\{\P^*_\mu\}$ is a family of elements in 
$A^*$ indexed by the vertices $\mu\in\G$, $\P^*_\varnothing=1$. We 
assume that these data obey the following  
condition which is a generalization of \tht{1.5}: 
$$ 
\P^*_\mu \P^*_{(1)}=a_n \P^*_\mu+  
\sum_{\la:\,\la\searrow\mu} \ka(\mu,\la)\P^*_\la\,, \qquad n=|\mu|, 
\tag 1.8 
$$ 
for any $\mu$, where $a_0=0,a_1,a_2,\dots$ is a sequence of numbers.  
 
\proclaim{Proposition 1.3} Under the above assumptions, let $\pi:A^*\to\C$ 
be a multiplicative linear functional, and let   
$$ 
s=\pi(\P^*_{(1)}), \quad t=-s=-\pi(\P^*_{(1)}).  
\tag 1.9 
$$  
Assume that  
$$ 
s\ne0,a_1,a_2,\dots\,, \quad \text{i.e.}, \quad  
t\ne0,-a_1,-a_2,\dots\,.   
\tag 1.10 
$$  
Then the function 
$$ 
\varphi(\mu)=\frac{\pi(\P^*_\mu)}{s(s-a_1)\dots(s-a_{n-1})} 
=\frac{(-1)^n\pi(\P^*_\mu)}{t(t+a_1)\dots(t+a_{n-1})}\,,   
\qquad n=|\mu|, 
\tag 1.11 
$$ 
is harmonic on $\G$. 
\endproclaim 
 
We agree that the denominator in \tht{1.11} equals 1 for 
$\mu=\varnothing$, so that $\varphi(\varnothing)=1$. 
 
\demo{Proof} Applying $\pi$ to the relation \tht{1.8} 
we get 
$$ 
\pi(\P^*_\mu)(s-a_n)=\sum_{\la:\, \la\searrow\mu}  
\ka(\mu,\la)\pi(\P^*_\la).  
$$ 
Dividing the both sides by $s(s-a_1)\dots(s-a_n)$ (which is possible 
thanks to \tht{1.10}) we get exactly the harmonicity relation \tht{1.1} for 
$\varphi$. \qed  
\enddemo 
 
A trivial example is $A^*=A$, $\P^*_\mu=\P_\mu$, $a_n\equiv0$. Then, by 
Theorem 1.2, $\varphi$ is extreme provided that it is nonnegative. As 
we aim to construct interesting examples of nonextreme harmonic 
functions, we shall deal either with an algebra $A^*$ distinct from 
$A$ or, for $A^*=A$, with a family $\{\P^*_\mu\}$ distinct from 
$\{\P_\mu\}$.  
 
In all the examples below, $A^*$ is a filtered 
algebra such that the associated graded algebra $\operatorname{gr}A^*$ is 
canonically isomorphic to $A$. Thus, with any element of $A^*$ of 
degree $\le n$ one can associate its {\it highest term\/} which is a 
homogeneous element of $A$ of degree $n$. In our examples, the 
highest term of $\P^*_\mu$ coincides with $\P_\mu$. Furthermore, 
the algebra $A^*$ can be interpreted, in a certain natural way, as an 
algebra of functions on the vertices of $\G$. Thus, for any $f\in 
A^*$ and $\la\in\G$, the value $f(\la)$ is well--defined. It turns out 
that the elements $\P^*_\mu$ can be characterized by the following  
 
\example{Interpolation Property} Given $\mu\in\G$,  
$\mu\ne\varnothing$, $\P^*_\mu$ is the only (up to a scalar factor) 
element of degree $|\mu|$ such that  
$\P^*_\mu(\la)=0$ for any $\la\ne\mu$ with $|\la|\le|\mu|$. 
\endexample 
 
The fact that the highest term of an element $\P^*_\mu$ defined in 
this way turns out to be proportional to $\P_\mu$ seems to be rather 
surprising. We normalize $\P^*_\mu$ in such a way that its highest 
term is exactly equal to $\P_\mu$.  
 
Next, it turns out that $\P^*_{(1)}(\mu)=|\mu|$. Then a simple formal 
argument shows that \tht{1.8} holds with $a_n=n$ for any $n=0,1,\dots$ . 
Moreover,  
$$ 
\frac{\dimG(\mu,\la)}{\dimG\la}=\frac{\P^*_\mu(\la)}{N(N-1)\dots(N-n+1)}\,, 
\qquad \mu\in\G_n, \quad \la\in\G_N, \quad n\le N. 
\tag 1.12 
$$ 
The argument is due to Okounkov \cite{Ok1}; it is also reproduced in 
\cite{OO}.  
 
{}From now on we shall assume that $a_n=n$. 
Then the denominator in the right--hand side of \tht{1.11} will be equal 
to $(t)_n=t(t+1)\cdots (t+n-1)$, and \tht{1.11} will take the form 
$$ 
\varphi(\mu)=\frac{(-1)^n\,\pi(\P^*_\mu)}{(t)_n}\,, 
\qquad t=-\pi(\P^*_{(1)}), \quad n=|\mu|.  
\tag 1.13 
$$ 
Similarly, the formula \tht{1.12} can be rewritten as follows 
$$ 
\frac{\dimG(\mu,\la)}{\dimG\la}=\frac{(-1)^n\, 
\P^*_\mu(\la)}{(-N)_n}\,,  
\qquad \mu\in\G_n, \quad \la\in\G_N, \quad n\le N.  
\tag 1.14 
$$ 
 
Note that, for any fixed $\la$, the left--hand side of \tht{1.14} satisfies 
the harmonicity relation \tht{1.1} provided that $n<N$: this easily follows 
from the very definition of the dimension function (for $n>N$ the 
denominator in the right--hand side vanishes). On the 
other hand, the expression in the right--hand side of \tht{1.14} is a 
particular case of that in the right--hand side of \tht{1.13}: here $\pi$ 
is the evaluation functional $\pi_\la: f\mapsto f(\la)$ and $t=-N$. 
This makes it possible to interpret the construction of Proposition 
1.3 as follows: we extrapolate the relation \tht{1.14} from the points 
$\la\in\G$, which we identify with the corresponding evaluation 
functionals $\pi_\la$, to abstract multiplicative functionals.  
 
A function $\varphi\in\Harm^+_1(\G)$ will be called {\it nondegenerate\/}  
if $\varphi(\mu)\ne0$ for all $\mu\in\G$; otherwise it will be called {\it 
degenerate}.

\head \S2. The Young graph \endhead 
 
The fundamental example of a graded graph $\G$ is the {\it Young graph\/} 
$\Y$ \cite{VK}, \cite{KV2}. By definition, the vertices of $\Y$ are the Young 
diagrams including the empty diagram $\varnothing$, the $n$-th floor 
$\Y_n$ consists of the diagrams with $n$ boxes, and $\mu\nearrow\la$ 
means that $\la$ is obtained from $\mu$ by adding a single box. The  
numbers $\ka(\mu,\la)$ are all equal to 1. In this section the 
symbols $\mu,\la$ are used to denote Young diagrams.  
 
The graph $\Y$ is multiplicative in the sense of the definition given 
in \S1: here the algebra $A$ 
is the algebra $\La$ of symmetric functions, the basis elements 
$\P_\mu$ are the Schur functions $s_\mu$, and the relation \tht{1.5} 
turns into a special case of the Pieri rule for the Schur functions, 
$$ 
s_\mu s_{(1)}=\sum_{\la:\,\la\searrow\mu} s_\la\,, 
\tag 2.1 
$$ 
which is equivalent (under the characteristic map, see \cite{M, I.7}) to 
the Young branching rule for irreducible characters of symmetric 
groups. For the Young graph, the expansion \tht{1.6}  
takes the form 
$$ 
s_{(1)}^n=\sum_{\la:\, |\la|=n} 
\dim\la\cdot s_\la,  
\tag 2.2 
$$ 
where $\dim\la=\dim_\Y\la$ is the number of standard Young tableaux of  
shape $\la$.  
 
Let 
$b=(i,j)$ be a box of $\mu$; here $i,j$ are the row number and  
the column number of $b$. Recall the definition of the {\it 
content\/}, the {\it arm--length\/} and the {\it  
leg--length\/} of $b$:  
$$ 
c(b)=j-i,\quad 
a(b)=\mu_i-j,\quad 
l(b)=\mu'_j-i,  
\tag 2.3 
$$ 
where $\mu'$ is the transposed diagram. 
 
\proclaim{Theorem 2.1} Let $z,z'$ be arbitrary complex numbers 
and $t=zz'$. Assume that $t\ne0,-1,-2,\dots$ . Then the 
following expression is a harmonic function on the Young graph: 
$$ 
\varphi_{zz'}(\mu)=\frac1{(t)_n} 
\prod_{b\in\mu}\frac{(z+c(b))(z'+c(b))}{a(b)+l(b)+1}\,,  
\quad n=|\mu|. 
\tag 2.4 
$$ 
 
The harmonic functions \tht{2.4} fit into the general scheme of 
Proposition 1.3 with the algebra $A^*$ and the family $\{\P^*_\mu\}$ 
as specified below.  
\endproclaim 
 
The first claim of the theorem (harmonicity of $\varphi_{zz'}$) 
follows from the computation of a spherical function in \cite{KOV}. 
Various direct combinatorial proofs for this claim were given by 
Kerov, Postnikov, and Borodin. Kerov's approach is explained in 
\cite{Ke5}; actually, in that paper a more general result is 
obtained, see Theorem 3.1 below. Postnikov's argument was not 
published. Borodin's argument is, perhaps, the most direct and 
elementary; it was given in the appendix to \cite{P.I}; 
actually, the present paper originated from our discussion of that 
argument.   
 
For the proof we need some preparations. First, we  specify the 
algebra $A^*$. 
 
Denote by $\La^*(n)$ the subalgebra in $\C[x_1,\dots,x_n]$ formed by 
the polynomials which are symmetric in `shifted' variables 
$x'_j=x_j-j$, $j=1,\dots,n$. Define the projection map 
$\Lt(n)\to\Lt(n-1)$ as the specialization $x_n=0$ and note that 
this projection preserves the filtration defined by ordinary degree of 
polynomials. Now we take the projective limit of $\La^*(n)$'s  
in the category of filtered algebras as $n\to\infty$. The result is a 
filtered algebra which is called the {\it algebra of shifted 
symmetric functions\/} and denoted by $\La^*$.  
 
The algebra $\La^*$ will be taken as the algebra $A^*$. As the 
elements $\P^*_\mu$ we shall take the {\it shifted Schur functions\/} 
$s^*_\mu$ as defined in \cite{OO}.  
 
By the definition of $\La^*$, each element $f\in\La^*$ can be 
evaluated at any sequence $x=(x_1,x_2, \dots)$ with finitely many 
nonzero terms. In particular, we can evaluate shifted symmetric 
functions at any $\la=(\la_1,\la_2,\dots)\in\Y$, which allows one to 
interpret $\La^*$ as a certain algebra of functions on the Young 
diagrams. This point of 
view was developed in \cite{KO}. The shifted Schur 
functions $s^*_\mu$ possess the Interpolation Property of \S1, see \cite{Ok1}, \cite{OO}.   
 
For the one--row shifted Schur functions there is a special notation: 
$h^*_m=s^*_{(m)}$. A useful tool is the following generating series for 
the $h^*$ functions: 
$$ 
H^*(u)=1+\sum_{m=1}^\infty \frac{h^*_m}{u(u-1)\dots(u-m+1)}\,. 
\tag 2.5 
$$ 
Here $u$ is a formal indeterminate and the series is viewed as an 
element of $\La^*[[\frac1u]]$. Since the elements $h^*_m$ are 
algebraically independent generators of $\La^*$, a multiplicative 
functional $\pi:\La^*\to\C$ can be uniquely defined by assigning to 
$H^*(u)$ an arbitrary formal power series in $\frac1u$ with constant 
term 1. We shall use this fact below.  
 
Note a useful formula 
$$ 
H^*(u)(x_1,x_2,\dots)=\prod_{i=1}^\infty  
\frac{u+i}{u+i-x_i}\,,  
\tag 2.6 
$$ 
see \cite{OO, Theorem 12.1}. Here, by definition, 
$$ 
H^*(u)(x_1,x_2,\dots)= 
1+\sum_{m=1}^\infty  
\frac{h^*_m(x_1,x_2,\dots)}{u(u-1)\dots(u-m+1)}\,.  
\tag 2.7 
$$ 
 
The equality \tht{2.6} can be understood as 
follows. We assume that only finitely many of $x_i$'s are distinct 
from zero. Then the left--hand side, which is the series \tht{2.7}, 
converges in a left half--plane $\Re u<\operatorname{const}\ll0$ and 
equals the right--hand side of \tht{2.6}.  
 
For an element $f$ of $\La$ or $\La^*$, we shall 
abbreviate  
$$ 
f(x_1,\dots,x_k)=f(x_1,\dots,x_k,0,0,\dots). 
$$ 
 
 Recall the 
combinatorial formula for the Schur functions: 
$$ 
s_\mu(x_1,\dots,x_k)=\sum_T \prod_{b\in\mu}x_{T(b)}\,,   
\tag 2.8 
$$ 
where $T$ ranges over the set of Young tableaux of shape $\mu$ 
with entries in $\{1,\dots,k\}$, see \cite{M, I.5}. It 
will be convenient for us to employ here the reverse tableaux (i.e., 
the entries $T(b)$ decrease from left to right along the rows and 
down the columns). Since $s_{\mu}$ is symmetric, (2.8) also holds if 
the sum in the right--hand side is taken over all reverse tableaux of 
shape $\mu$ with entries in $\{1,\dots,k\}$.   
 
 We shall need a similar formula for the shifted 
Schur functions: 
$$ 
s^*_\mu(x_1,\dots,x_k)=\sum_T \prod_{b\in\mu}(x_{T(b)}-c(b))\,,   
\tag 2.9 
$$ 
where $T$ ranges over reverse tableaux of shape $\mu$ with entries in 
$\{1,\dots,k\}$, see \cite{OO, Theorem 11.1}.  
 
\proclaim{Proposition 2.2} Let $k=1,2,\dots$ and $z'\in\C$. The 
following specialization formula holds 
$$ 
s^*_\mu(\underbrace{-z',\dots,-z'}_k)= 
(-1)^n\prod_{b\in\mu} 
\frac{(k+c(b))(z'+c(b))}{a(b)+l(b)+1}\,, 
\quad n=|\mu|.  
\tag 2.10 
$$ 
\endproclaim 
 
\demo{Proof} Compare the combinatorial formulas \tht{2.8} and \tht{2.9}. 
If $x_1=\dots=x_k=-z'$ then the product in \tht{2.9} does not 
depend on $T$ and is equal to 
$$ 
\prod_{b\in\mu}(-z'-c(b)) 
=(-1)^n\prod_{b\in\mu}(z'+c(b)). 
$$ 
It follows that  
$$ 
s^*_\mu(\underbrace{-z',\dots,-z'}_k) 
=(-1)^n\prod_{b\in\mu}(z'+c(b))\cdot  
s_\mu(\underbrace{1,\dots,1}_k). 
$$ 
 
Now we apply the well--known specialization formula  
$$ 
s_\mu(\underbrace{1,\dots,1}_k) 
=\prod_{b\in\mu}\frac{k+c(b)}{a(b)+l(b)+1}\,, 
$$ 
see \cite{M, I.3, Ex.4}, which implies \tht{2.10}. \qed 
\enddemo 

The argument used in the proof is borrowed from Okounkov's paper 
\cite{Ok4}, the derivation of formula (1.9); see also Proposition 3.2 below.

\proclaim{Corollary 2.3} For any $z,z'\in\C$, the linear functional 
$\pi_{zz'}:\La^*\to\C$ given by  
$$ 
\pi_{zz'}(s^*_\mu)= 
(-1)^n\prod_{b\in\mu} 
\frac{(z+c(b))(z'+c(b))}{a(b)+l(b)+1}, 
\quad n=|\mu|,  
\tag 2.11 
$$ 
is multiplicative. 
\endproclaim 
 
\demo{Proof} Indeed, this expression depends polynomially on $z$. So, 
it suffices to prove the multiplicativity of $\pi_{zz'}$ in  
the case $z=k$, where $k=1,2,\dots$ . By Proposition 2.2, in this 
case our functional is the evaluation at the point 
$$
x=(\underbrace{-z',\dots,-z'}_k,0,0,\dots).
$$
 Consequently, the functional is multiplicative. \qed  
\enddemo 
 
\demo{Proof of Theorem 2.1} We apply Proposition 1.3 by taking 
$A^*=\La^*$ and $\P^*_\mu=s^*_\mu$. The Pieri--type formula for 
$s^*$--functions ([OO, Theorem 9.1]) shows that the relation 
\tht{1.8}  
holds with the sequence $a_n=n$. We take as $\pi$ the multiplicative 
functional $\pi_{zz'}$ afforded by Corollary 2.3. It follows from 
\tht{2.11} that  
$$ 
-\pi_{zz'}(s^*_{(1)})=zz'=t, 
$$ 
so that we may substitute $t$ into \tht{1.11}. Finally, the condition 
\tht{1.10} is just the assumption on $t$ given in Theorem 2.1. Thus, the 
expression \tht{2.4} is a special case of \tht{1.11}, which concludes the 
proof. \qed 
\enddemo 
 
\example{Remark 2.4} In terms of the generating series \tht{2.5} for the 
$h^*$ functions, the multiplicative functional $\pi_{zz'}$ can be 
described as follows: 
$$ 
\gathered 
\pi_{zz'}(H^*(u))=1+\sum_{m=1}^\infty\frac{(z)_m(z')_m}{(-u)_m\,m!}\\ 
={}_2F_1(z,z';-u;1) 
=\frac{\Ga(-u)\Ga(-u-z-z')}{\Ga(-u-z)\Ga(-u-z')}\,,   
\endgathered 
\tag 2.12 
$$ 
where we assume $\Re u\ll0$; the last equality follows from Gauss' 
summation formula.  
\endexample 
 
\proclaim{Proposition 2.5} The function $\varphi_{zz'}$ afforded by 
Theorem 2.1 is a nondegenerate function from $\Harm^+_1(\Y)$ if and 
only if the parameters satisfy one of the following two conditions: 
 
$\bullet$ either $z'=\bar z$ where $z\in\C\setminus\Z$, 
 
$\bullet$ or $z,z'$ are real and there exists $m\in\Z$ such that 
$m<z,z'<m+1$.  
\endproclaim 

The proof is easy, see \cite{P.I}.

Let us explain the significance of the set $\Harm^+_1(\G)$ for the 
Young graph. Let $S(\infty)$ be the infinite symmetric group, which 
is defined as the inductive limit $\varinjlim S(n)$ of the finite 
symmetric groups. In other words, $S(\infty)$ consists of the finite 
permutations of the set $\{1,2,\dots\}$. There is a natural bijective 
correspondence $\varphi\leftrightarrow\chi$ between  functions 
$\varphi\in\Harm(\Y)$ and central functions $\chi$ on 
$S(\infty)$. Specifically, given $\chi$, we define the values of 
$\varphi$ on $\Y_n$ from the expansion of the central function 
$\chi\downarrow S(n)$ on the group $S(n)$ into a linear combination 
of the irreducible characters $\chi^\la$, 
$$ 
\chi\downarrow S(n)=\sum_{\la:\, |\la|=n}\varphi(\la)\chi^\la\,, 
\qquad n=1,2,\dots\,.  
\tag 2.13 
$$ 
The harmonicity of the function $\varphi$ follows from the  Young 
branching rule for the irreducible characters \cite{JK}, \cite{OV}: 
$$ 
\chi^\la\downarrow S(n-1)=\sum_{\mu:\, \mu\nearrow\la}\chi^\mu\,, 
\quad n=|\la|.  
\tag 2.14 
$$  
Under the bijection $\varphi\leftrightarrow\chi$, the nonnegativity of 
$\varphi$  means that the function $\chi$ 
is positive definite on the group $S(\infty)$, and the 
normalization $\varphi(\varnothing)=1$ means that $\chi(e)=1$, where 
$e\in S(\infty)$ is the trivial permutation. Thus, the elements of 
$\Harm^+_1(\Y)$ correspond to the central, positive definite, 
normalized functions on $S(\infty)$. Such functions form a convex 
set, its extreme points are called the {\it characters} of the group 
$S(\infty)$ (in the sense of von Neumann). According to \S1, we shall 
denote the set of the characters by $\Om(\Y)$.   
 
Thoma proved [T] that the characters of $S(\infty)$ can be 
parametrized by the points of an infinite--dimensional simplex:  
$$  
\Om(\Y)=\{\al_1\ge\al_2\ge\dots\ge 0,\ \beta_1\ge\beta_2\ge\dots\ge 
0\,|\, \sum_i(\alpha_i+\beta_i)\le 1\}  
\tag 2.15 
$$   
called the {\it Thoma simplex}. It is equipped with the weakest topology in which the coordinates $\alpha_i$'s and $\beta_i$'s are continuous 
functions.  
  
Via Gelfand--Naimark--Segal construction, characters generate {\it 
finite factor representations} of the group $S(\infty)$. They also 
correspond to {\it irreducible unitary spherical representations} of 
the Gelfand pair $(G,K)$ where $G$ is the ``bisymmetric group'' 
$S(\infty)\times S(\infty)$, and $K$ is the diagonal subgroup of $G$ 
\cite{Ol}, \cite{Ok3}.  
 
The Poisson kernel $K(\mu,\om)$ (see Theorem 1.1) for the Young graph 
is given by the image of the Schur function $s_\mu$ under a certain 
specialization of $\Lambda$ depending on $\omega$. Namely, for 
$\omega=(\alpha,\beta)\in\Omega(\Y)$ we specialize the power sums as 
follows  
$$ 
(x_1^k+x_2^k+\ldots)\mapsto\cases 1,& k=1,\\  
\sum_{i=1}^\infty\alpha_i^k+(-1)^{k-1} 
\sum_{i=1}^{\infty}\beta_i^k,& k\ge 2, 
\endcases 
\tag 2.16 
$$  
see \cite{VK}, \cite{KOO}. 
 
According to the general theory of \S1, the functions 
$\varphi_{zz'}\in\Harm_1^+(\Y)$ constructed above give rise to a 
family of probability measures on $\Omega(\Y)$. These measures were 
thoroughly studied in \cite{P.I--P.V}, \cite{BO1}, \cite{BO2}.

\head \S3 The Jack graph \endhead 
 
Fix a positive number $\th$. Let $P_\mu$ be the Jack symmetric 
function with parameter $\th$ and index $\mu$ (see \cite{M, VI.10}; note 
that Macdonald uses $\alpha=\th^{-1}$ as the parameter). The simplest 
case of Pieri's formula for the Jack functions reads as follows: 
$$ 
P_\mu P_{(1)}=\sum_{\la:\, \la\searrow\mu}\dimth(\mu,\la)P_\la\,,   
\tag 3.1 
$$ 
where $\dimth(\mu,\la)$ are certain positive numbers, 
$$ 
\dimth(\mu,\la) = \prod_b \frac 
{\big(a(b) + (l(b)+2)\th\big)\big(a(b) + 1 + l(b)\th\big)} 
{\big(a(b)+(l(b)+1)\th\big)\big(a(b)+1+(l(b)+1)\th\big)}\,. 
\tag 3.2  
$$ 
Here $b$ ranges over all boxes in the $j$th column of the diagram 
$\mu$, provided that the new box $\la\setminus\mu$ belongs to the 
$j$th column of $\la$, see \cite{M, VI.10, VI.6}.  
 
The {\it Jack graph\/} $\J$ is the multiplicative graph associated with 
the algebra $A=\La$ of symmetric functions and its basis formed by 
the Jack symmetric functions. I.e., this is the Young graph with the 
formal edge multiplicities $\ka(\mu,\la)=\dimth(\mu,\la)$. When 
$\th=1$, the Jack functions turn into the Schur functions, all  
formal edge multiplicities are equal to 1, so that the Jack graph 
turns into the ordinary Young graph. 
 
We take as $A^*$ the algebra $\Lt$ of shifted symmetric 
functions with parameter $\th$. It is defined as the projective limit 
of filtered algebras $\Lt(n)$, where, in turn, $\Lt(n)$ is formed by 
polynomials in $x_1,\dots,x_n$ which are symmetric with respect to 
the new variables $x'_j=x_j-\th j$. The projection 
$\Lt(n)\to\Lt(n-1)$, as in the case of the Young graph, is given by 
the specialization $x_n=0$.  
 
For any $\mu$ there exists a unique element $P^*_\mu\in\Lt$ of degree 
$|\mu|$, with highest term $P_\mu$ and with the Interpolation 
Property; it is called the {\it shifted Jack function.\/} Our 
reference about these functions is \cite{OO2}, \cite{Ok2}. 
 
The  
relation \tht{1.8} for shifted Jack functions has the form  
$$ 
P^*_\mu P^*_{(1)} 
=nP^*_\mu+\sum_{\la:\, \la\searrow\mu}\dimth(\mu,\la)P^*_\la\,, 
\qquad n=|\mu|.   
\tag 3.3 
$$ 
 
Given a box $b=(i,j)$ of a Young diagram, we denote its $\th$-{\it 
content\/} as $c_\th(b)=(j-1)-\th(i-1)$. When $\th=1$, this turns into 
the ordinary content. 
 
\proclaim{Theorem 3.1} Let $z,z'$ be arbitrary complex numbers 
and $t=\th^{-1}zz'$. Assume that $t\ne0,-1,-2,\dots$ . Then the 
following expression is a harmonic function on the Jack graph 
$\J${\rm :} 
$$ 
\varphi_{zz'}(\mu)=\frac1{(t)_n} 
\prod_{b\in\mu}\frac{(z+c_\th(b))(z'+c_\th(b))}{a(b)+\th 
l(b)+\th}\,, \quad n=|\mu|. 
\tag 3.4 
$$ 
 
The functions \tht{3.4} fit into the general scheme of Proposition 1.3 
with $A^*=\La^\th$ and $\P^*_\mu=P^*_\mu$. 
\endproclaim 
 
This result generalizes Theorem 2.1. The first claim is due to Kerov 
\cite{Ke5}. Our proof of Theorem 3.1 is very similar to that of Theorem 
2.1, so we shall only indicate necessary modifications.  
 
The analogs of the combinatorial formulas \tht{2.8} and \tht{2.9} are as follows:  
$$ 
\gather 
P_\mu(x_1,\dots,x_k)=\sum_T \psi_T(\th) 
\prod_{b\in\mu}x_{T(b)}\,,  
\tag 3.5 \\ 
P^*_\mu(x_1,\dots,x_k)=\sum_T \psi_T(\th) 
\prod_{b\in\mu}(x_{T(b)}-c_\th(b))\,,  
\tag 3.6 
\endgather 
$$ 
where, again, the summation is taken over the reverse Young tableaux of 
shape $\mu$, and $\psi(\th)$ are certain numeric factors. We do not 
need their exact values, the point is that they are the same in both 
formulas, see \cite{Ok2}, \cite{Ok4}.  
 
\proclaim{Proposition 3.2} Let $k=1,2,\dots$ and $z'\in\C$. The 
following specialization formula holds 
$$ 
P^*_\mu(\underbrace{-z',\dots,-z'}_k)= 
(-1)^n\prod_{b\in\mu} 
\frac{(k+c_\th(b))(z'+c_\th(b))}{a(b)+\th l(b)+\th}\,, 
\quad n=|\mu|.  
\tag 3.7 
$$ 
\endproclaim 
 
\demo{Proof} The argument is exactly similar to that for Proposition 
2.2. We employ the combinatorial formulas \tht{3.5}, \tht{3.6} and the 
well--known specialization formula for the Jack symmetric functions: 
$$ 
P_\mu(\underbrace{1,\dots,1}_k) 
=\prod_{b\in\mu}\frac{\th k+c_\th(b)}{a(b)+\th l(b)+\th}\,, 
$$ 
see \cite{M, VI, (10.20)}. \qed 
\enddemo 

A more general formula appeared in \cite{Ok4, (1.9)}.

Proposition 3.2 immediately leads to the following claim. 
 
\proclaim{Corollary 3.3} For any $z,z'\in\C$, the linear functional 
$\pi_{zz'}:\Lt\to\C$ given by  
$$ 
\pi_{zz'}(P^*_\mu)= 
(-1)^n\prod_{b\in\mu} 
\frac{(z+c_\th(b))(z'+c_\th(b))}{a(b)+\th l(b)+\th}, 
\quad n=|\mu|,  
\tag 3.8 
$$ 
is multiplicative. 
\endproclaim 
 
\demo{Proof} The argument is just the same as for Corollary 2.3. \qed 
\enddemo 
 
\demo{Proof of Theorem 3.1} Exactly the same as for Theorem 2.1. \qed 
\enddemo 
 
Extreme points of the convex set $\Harm_1^+(\J)$ of nonnegative 
normalized harmonic functions on the Jack graph, as in the case of 
the Young graph ($\theta=1$), can be parametrized by points of the 
Thoma simplex (2.15), see \cite {KOO}.   
 
The Poisson kernel $K(\mu,\omega)$ is defined as the image of the 
Jack function $P_\mu$ under the specialization of $\Lambda$ which 
sends power sums to the following expressions (cf. (2.16)):  
 $$ 
(x_1^k+x_2^k+\ldots)\mapsto\cases 1,& k=1,\\  
\sum_{i=1}^\infty\alpha_i^k+(-\theta)^{k-1} 
\sum_{i=1}^{\infty}\beta_i^k,& k\ge 2, 
\endcases 
\tag 3.9 
$$  
see \cite{KOO}. 
 
As was mentioned in \S2, the set $\Harm_1^+(\Y)=\Harm_1^+(\JJ(1))$ has a representation theoretic meaning. There is one more value of $\theta$, namely $\theta=1/2$, when harmonic functions on the Jack graph can be related to representations. We shall briefly explain this connection.

Let $G$ be the group of finite permutations of the set $\{\pm 1,\pm 2,\dots\}$ and $K$ be its subgroup consisting of the permutations which commute with the involution $i\mapsto -i$. The group $G$ is just another realization of the infinite symmetric group $S(\infty)$, and the group $K$ can be considered as an infinite version of the hyperoctahedral groups.
Note that $(G,K)$ is a Gelfand pair \cite{Ol}. 

It turns out that there exists a natural one--to--one correspondence between  $\Harm_1^+(\JJ(1/2))$ and the set of positive definite, $K$--biinvariant functions on $G$ normalized at the unity (this correspondence is based on classical facts explained in \cite{M, VII.2}). In particular, extreme functions from $\Harm_1^+(\JJ(1/2))$ correspond to the spherical functions of irreducible unitary spherical representations of $(G,K)$. For more details about these representations, see \cite{Ol}, \cite{Ok3}.

\head 4. The Kingman graph \endhead  
 
Let 
$$ 
\mu=(\mu_1\ge\dots\ge\mu_l>0)=(1^{r_1(\mu)}2^{r_2(\mu)}\dots) 
$$ 
denote an arbitrary partition also viewed as a Young diagram.  
 
In this section we are dealing with the monomial symmetric functions 
$m_\mu$ \cite{M, I.2}. They form a basis of the algebra  
$\La$ and obey the relation: 
$$ 
m_\mu m_{(1)}=\sum_{\la:\, \la\searrow\mu} \ka_0(\mu,\la)m_\la\,, 
\tag  4.1 
$$ 
where the positive integers $\ka_0(\mu,\la)$ are defined as follows: 
if $k$ stands for the length of the 
row in $\la$ containing the box $\la\setminus\mu$ then  
$\ka_0(\mu,\la)=r_k(\la)$.  
 
The {\it Kingman graph\/} $\K$ is the multiplicative graph associated 
with the algebra $\La$ and its basis $\{m_\mu\}$ \cite{Ke1}. I.e., this is the 
Young graph with the formal edge multiplicities $\ka_0(\mu,\la)$. Since 
the numbers $\ka_0(\mu,\la)$ are integers, one can regard $\K$ as a 
graph with multiple edges.  
 
Next, introduce the {\it factorial monomial symmetric functions\/} 
$m^*_\mu$, which are also elements of $\La$. By definition 
\cite{Ke1}, $m^*_\mu$ is the sum of all 
{\it distinct\/} expressions obtained from 
$$ 
\prod_{i=1}^l x_i(x_i-1)\dots(x_i-\mu_i+1) 
$$ 
by permutations of the variables $x_1,x_2,\dots$ . Thus, the definition of $m^*_\mu$ is similar to that of $m_\mu$, the only difference 
is that the ordinary powers $x^m$ are replaced by the 
falling factorial powers $x(x-1)\dots(x-m+1)$. 
 
The function $m^*_\mu$ can be characterized as the only symmetric function 
with the highest term $m_\mu$ and such that 
$m^*_\mu(\la_1,\la_2,\dots)=0$ for any diagram $\la\ne\mu$, 
$|\la|\le|\mu|$. Thus, $m^*_\mu$ possesses the Interpolation Property 
of \S1.  
 
One can directly verify that 
$$ 
m^*_\mu m^*_{(1)}=n m^*_\mu+\sum_{\la:\, \la\searrow\mu} 
\ka_0(\mu,\la)m^*_\la, \qquad n=|\mu|   
\tag 4.2 
$$ 
(this also follows from the Interpolation Property).  
 
\proclaim{Theorem 4.1} Let $t,\al$ be complex parameters, 
$t\ne-1,-2,\dots$\,.  
Then the function 
$$ 
\varphi_{t,\al}(\mu)= 
\frac{(\mu_1-1)!\dots(\mu_l-1)!}{r_1(\mu)!r_2(\mu)!\dots}\cdot 
\frac{t(t+\al)\dots(t+(l-1)\al)}{(t)_n}\cdot 
\prod\Sb b=(i,j)\in\mu\\j\ge2\endSb\left(1-\frac{\al}{j-1}\right)\,, \tag 4.3 
$$ 
where $n=|\mu|$, is harmonic on the graph $\K$. Here $l$ is the length (number of nonzero parts) of $\mu$.  
 
The functions \tht{4.3} fit into the general scheme of Proposition 1.3 with 
$A^*=\La$ and $\P^*_\mu=m^*_\mu$. 
\endproclaim 
 
As is explained below, the first claim is equivalent to a result of 
Pitman \cite{Pi}.  
 
\demo{Proof} According to Proposition 1.3, 
it suffices to check that  
there exists a multiplicative functional $\pi_{t,\al}: \La\to\C$ such 
that  
$$ 
\gathered 
\pi_{t,\al}(m^*_\mu)=(-1)^n\, 
\frac{(\mu_1-1)!\dots(\mu_l-1)!}{r_1(\mu)!r_2(\mu)!\dots}\\ \cdot\,
t(t+\al)\dots(t+(l-1)\al)\cdot 
\prod\Sb b=(i,j)\in\mu\\j\ge2\endSb\left(1-\frac{\al}{j-1}\right)\,.  \endaligned 
\tag 4.4 
$$ 
 
As the functions $m^*_\mu$ form a basis in $\La$, we can define a 
linear functional $\pi_{t,\al}:\La\to\C$ by the formula \tht{4.4}. We 
claim that it is multiplicative if $t=-k\al$, where $k=1,2,\dots$\,. 
To see this we shall prove that $\pi_{-k\al,\al}$ coincides with the 
specialization at the point $(\underbrace{\al, \dots,\al}_k)$.  
 
Indeed, from the definition of $m^*_\mu$ it follows that   
$$ 
m^*_\mu(\underbrace{\al,\dots,\al}_k) 
=\frac{k(k-1)\dots(k-l+1)}{r_1(\mu)!r_2(\mu)!\dots}\, 
\prod_{i=1}^l \al(\al-1)\dots(\al-\mu_i+1), 
\tag 4.5    
$$ 
and a direct verification shows that this expression coincides with 
$\pi_{-k\al,\al}(m^*_\mu)$.  
 
Finally, as the right--hand side of \tht{4.4} depends on the 
parameters $t,\al$ polynomially, $\pi_{t,\al}$ is multiplicative for 
all values of the parameters. \qed  
\enddemo 
 
\proclaim{Proposition 4.2} The function $\varphi_{t,\al}$ afforded by 
Theorem 4.1 is a nondegenerate function from $\Harm^+_1(\K)$ if and 
only if the parameters $t,\al$ are real and satisfy the inequalities 
$0\le\al<1$, $t>-\al$.  
\endproclaim

The proof is straightforward.  
 
There is a bijective correspondence between the functions 
$\varphi\in\Harm^+_1(\K)$ and the {\it partition structures\/} in the 
sense of Kingman \cite{Ki1}, \cite{Ki2}. According to Kingman, a partition 
structure is a sequence $M=(M_n)$ of probability distributions on 
the partitions of $n$, $n=1,2,\dots$, such that for each $n$, $M_n$ 
and $M_{n+1}$ are connected by a certain consistency relation. These 
sequences are nothing else than the sequences $(M_n)$ as defined in 
\S1, see (1.4). Thus, the passage from a harmonic function 
$\varphi\in\Harm^+_1(\K)$ to the corresponding partition structure  
is given by the formula   
$$ 
M_n(\mu)=\varphi(\mu)\cdot\dim_\K(\mu), \qquad |\mu|=n, 
$$ 
where 
$$ 
\dim_\K(\mu)=\frac{n!}{\mu_1!\dots\mu_l!}  
\tag 4.6 
$$ 
(the latter formula readily follows from the general relation 
\tht{1.12} if we substitute $\lambda=\mu$ and $\Cal 
P^*_\mu=m^*_\mu$).   
 
Under this correspondence, the functions $\varphi_{t,0}$ with $t>0$ turn 
into Ewens' partition structures \cite{Ew}. More general 
functions $\varphi_{t,\al}$ with the restrictions $t>-\al$, $0\le\al<1$, 
of Proposition 4.2 correspond to Pitman's two--parameter 
generalization of Ewens' partition structures  
\cite{Pi}, \cite{Ke4}. Note that the harmonic functions $\varphi_{-k\al,\al}$ with 
$k=1,2,\dots$, which appear in the proof of Theorem 4.1, are nonnegative and 
{\it degenerate\/} provided that $\al<0$; the significance of the 
corresponding partition structures is explained in the introduction 
to \cite{Pi}.  
 
Here is yet another interpretation of the harmonic functions $\varphi_{t,\al}$. 
 
For $n\ge2$ there exists a unique map $S(n)\to S(n-1)$ 
which commutes with the two--sided action of the smaller group 
$S(n-1)$. This map, called the {\it canonical projection,\/} can be 
defined as follows: if $s,s_1,s_2\in S(n-1)$ and $(n-1,n)$ stands for 
the elementary transposition of ``$n-1$'' and ``$n$'', then 
$s\mapsto s$ and $s_1\cdot(n-1,n)\cdot s_2\mapsto s_1s_2$. In other 
words, the canonical projection is defined by removing ``$n$'' from 
the cycle $\dots i\to n\to j\dots$ containing it, see \cite{KOV}.  
 
The projective limit $\X=\varprojlim S(n)$ taken with respect to the 
canonical projections is a compact topological space. Its elements are 
called {\it virtual permutations.\/} There is a natural 
embedding $S(\infty)\to\X$ whose image is dense, so that $X$ is a 
certain compactification of the discrete set $S(\infty)$. The 
two--sided action of $S(\infty)$ on itself can be extended to $\X$, 
which makes $\X$ a $S(\infty)\times S(\infty)$-space. This 
construction and its meaning for the representation theory of the 
group $S(\infty)$ is discussed in \cite{KOV}.  
 
A probability measure $\M$ on $\X$ is called {\it central\/} if it is 
invariant under the action of the diagonal subgroup in  
$S(\infty)\times S(\infty)$ (that action extends the action of 
$S(\infty)$ on itself by conjugations).  
There is a natural bijective correspondence 
$\varphi\leftrightarrow \M$ between the elements  
$\varphi\in\Harm^+_1(\K)$ and the central measures $\M$ on $\X$. It is 
specified as follows: for each $n=1,2,\dots$, the image of $\M$ under 
the composite projection $\X\to S(n)\to\K_n$ coincides with $M_n$; here 
the second arrow $S(n)\to\K_n$ assigns to a permutation its cycle 
structure which is identified with a partition, and $(M_n)$ is the 
partition structure corresponding to $\varphi$.  
 
Let us denote by $\M_{t,\al}$ the measures corresponding to the 
harmonic functions $\varphi_{t,\al}$, where the parameters satisfy the 
conditions of Proposition 4.2. The measure $\M_{1,0}$ is invariant with 
respect to $S(\infty)\times S(\infty)$ and it is the only probability 
measure with this property. The measures $\M_{t,0}$ with $t>0$ were 
employed in \cite{KOV} for a geometric construction of generalized regular representations of the group $S(\infty)\times S(\infty)$ which are closely related to the harmonic functions $\varphi_{zz'}$ defined in \S2. Note that all the measures $\M_{t,\al}$ are quasi--invariant under the action of this group, see \cite{KOV}, \cite{Ke4}. 
 
Extreme points of the convex set $\Harm^+_1(\K)$ can be parametrized 
by nonincreasing sequences of nonnegative numbers with sum less 
or equal to 1 (cf. \tht{2.15}):  
$$ 
\Om(\K)=\{\al_1\ge\al_2\ge\ldots\,|\, \sum_{i}\al_i\le 1\}, 
\tag 4.7 
$$ 
and the Poisson kernel $K(\mu,\omega)$ is given in this case 
by {\it extended monomial symmetric functions}:  
$$ 
K(\mu,\al)=K(1^{r_1}2^{r_2}\dots,\al)= 
\sum_{k=0}^{r_1}\frac{(1-\sum_i\al_i)^k}{k!}\, 
m_{(1^{r_1-k}2^{r_2}\dots)}(\al_1,\al_2,\dots), 
\tag 4.8 
$$ 
see \cite{Ki2}, \cite{Ke4}.
  
Let us also note that the Kingman graph 
may be viewed as the degeneration of the Jack graph $\J$ as 
$\th\to0$. Indeed, according to the definition of the Jack functions,  
their expansion in the monomial functions has the form 
$$ 
P_\mu=m_\mu\,+\,\text{lower terms} 
$$  
relative to the dominance order on partitions \cite{M, VI, (10.13)}. It is well 
known that in the limit $\th\to 0$ the coefficients of all the lower 
terms vanish. In this sense, the Jack functions $P_\mu$ 
degenerate to the monomial functions $m_\mu$ as $\th\to0$. This 
implies, in particular, the limit relation  
$$ 
\ka_0(\mu,\la)=\lim_{\th\to0}\ka_\th(\mu,\la), 
$$  
which can also be checked directly from \tht{3.2}.

\head \S5. The Schur graph \endhead 
 
Recall that a partition is said to be {\it strict\/} if its nonzero 
parts are pairwise distinct. In this section, the symbols $\mu$ 
and $\la$ always mean strict partitions. Using the standard 
correspondence between partitions and Young diagrams we introduce the 
relation $\mu\nearrow\la$ as before, see \S2. Then the {\it Schur graph\/} 
$\SS$ is defined as follows: the vertices of the  
$n$th floor $\SS_n$ are the strict partitions of $n$, and the 
edges are the couples $\mu\nearrow\la$. By definition, the empty 
partition $\varnothing$ is included into the set of the strict 
partitions. All the edge multiplicities are equal to 1. Thus, the 
Schur graph is a subgraph of the Young graph.  
 
Let $\Ga$ denote the subalgebra in $\La$ generated by the odd power 
sums $p_1=\sum_i x_i$, $p_3=\sum_i x_i^3,\dots$\,. Equivalently, $\Ga$ 
consists of those symmetric functions $f(x_1,x_2,\dots)$ which 
satisfy the following cancellation condition: for any $i\ne j$, the 
result of specializing $x_i=y$, $x_j=-y$ in $f$ does not depend on 
$y$ (see \cite{Pr}).  
 
In this section, the symbol $P_\mu$ stands for the Schur $P$ function 
indexed by a strict partition $\mu$. The Schur $P$ functions form a 
homogeneous basis of $\Ga$, $\deg P_\mu=|\mu|$. They obey the 
following Pieri--type rule: 
$$ 
P_\mu P_{(1)}=\sum_{\la:\,\la\searrow\mu} P_\la\,,  
\tag 5.1 
$$ 
see \cite{M, III.8}. Thus, $\SS$ is a multiplicative graph with $A=\Ga$ and 
$\{\P_\mu\}=\{P_\mu\}$.   
 
Note that 
$$ 
\dim_\SS\mu=\frac{n!}{\prod_{i=1}^{\ell(\mu)}\mu_i!} 
\prod_{1\le i<j\le\ell(\mu)}\frac{\mu_i-\mu_j}{\mu_i+\mu_j}\,.  
\tag 5.2  
$$ 
Here and below $\ell(\mu)$ denotes the length of $\mu$. 
 
We shall employ the {\it factorial Schur $P$ functions\/} $P^*_\mu$, 
see \cite{I1}. These are inhomogeneous elements of $\Ga$ with the 
Interpolation Property; the highest term of $P^*_\mu$ coincides with 
$P_\mu$. According to the general formalism, the $P^*$ functions 
satisfy the relation  
$$ 
P^*_\mu P^*_{(1)}=n P^*_\mu+  
\sum_{\la:\,\la\searrow\mu} P^*_\la\,, \qquad n=|\mu|, 
\tag 5.3 
$$ 
which has the same form as for the shifted Schur functions (see 
[OO, Theorem 9.1]), except that now $\mu$ and $\la$ are not arbitrary 
but strict partitions.  
 
There is a convenient generating series for the one--row $P^*$ 
functions,  
$$ 
F^*(u)=1+\sum_{m=1}^\infty \frac{P^*_{(m)}}{u(u-1)\dots(u-m+1)}\,, 
\tag 5.4 
$$ 
whose evaluation at a point  $x=(x_1,x_2,\dots)$ has the form 
$$ 
F^*(u)(x)=\prod_{i=1}^\infty \frac{u+1+x_i}{u+1-x_i}\,,  
\tag 5.5 
$$ 
see \cite{I2}. The formulas \tht{5.4}, \tht{5.5} should be understood in the same way as the formulas \tht{2.6}, \tht{2.7}. 
 
\proclaim{Theorem 5.1} Let $t>0$. The 
following expression is a positive harmonic function on the Schur graph: $$ 
\varphi_t(\mu)=\frac1{(t)_n}\, 
\frac{\prod_{(i,j)\in\mu} (2t+(j-1)j)} 
{2^{\ell(\mu)}\,\prod_{i=1}^{\ell(\mu)}\mu_i!} \, 
\prod_{1\le i<j\le\ell(\mu)}\frac{\mu_i-\mu_j}{\mu_i+\mu_j}\,,  
\qquad n=|\mu|.  
\tag 5.6 
$$ 
 
The harmonic functions \tht{5.6} fit into the general scheme of 
Proposition 1.3 with $A^*=\Ga$ and $\P^*_\mu=P^*_\mu$.  
\endproclaim 
 
The first claim was established in \cite{B1}.  
 
Let $\widetilde\mu$ denote the {\it shifted\/} Young diagram 
corresponding to $\mu$ (see \cite{M, III.8}). Define $t_1,t_2$ 
from the conditions 
$$ 
t_1+t_2=1,\qquad t_1t_2=2t.   
\tag 5.7 
$$ 
Then (cf. \tht{2.4}) 
$$ 
\prod_{(i,j)\in\mu} (2t+(j-1)j)=\prod_{b\in\widetilde\mu} 
(t_1+c(b))(t_2+c(b))\,  
\tag 5.8 
$$ 
so that the parameters $t_1,t_2$ are to a certain extent similar to 
$z,z'$.  
 
\demo{Sketch of proof} As usual, we shall employ Proposition 1.3.  
 
Consider the linear functional $\pi_t:\Ga \to\C$ 
defined by 
$$ 
\pi_t(P^*_\mu)=(-1)^n\, 
\frac{\prod_{(i,j)\in\mu} (2t+(j-1)j)} 
{2^{\ell(\mu)}\,\prod_{i=1}^{\ell(\mu)}\mu_i!} \, 
\prod_{1\le i<j\le\ell(\mu)}\frac{\mu_i-\mu_j}{\mu_i+\mu_j}\,,  
\qquad n=|\mu|.  
\tag 5.9 
$$ 
Here $t\in\C$ is arbitrary. Note that $\pi_t(P^*_{(1)})=-t$, as 
prescribed by Proposition 1.3. According to the general formalism, it suffices to prove 
that $\pi_t$ is multiplicative.  
 
Since the expression \tht{5.9} depends on $t$ polynomially, it suffices 
to prove the claim for a countable number of different values of $t$. 
We shall assume that  
$$ 
t=\frac{k(1-k)}2\,, \qquad k=1,2,\dots\, .  
\tag 5.10 
$$ 
That is, in the notation of \tht{5.7},  
$$ 
t_1=k, \quad t_2=1-k, \qquad k=1,2,\dots\, .  
\tag 5.11 
$$ 

To prove the multiplicativity property for the values \tht{5.10} we 
shall show that  
$$ 
\pi_{k(1-k)/2}=\text{  
evaluation at the staircase 
diagram } (k,k-1,\dots,1). 
\tag 5.12 
$$ 
 
This claim is an analog of Propositions 2.2, 3.2, and formula (4.5). 
It does not seem to have appeared in the literature before, so we give here a 
sketch of the proof.  
 
We shall consecutively check 
\tht{5.12} on the one--row $P^*$ functions, next on the two--row $P^*$  
functions, and finally on arbitrary $P^*$ functions. 
 
First, consider the generating series \tht{5.4} for the one--row functions. By 
\tht{5.5}, its evaluation at the $k$th staircase diagram is as follows: 
$$ 
F^*(u)(k,k-1,\dots,1)=\prod_{i=1}^k\frac{u+1+i}{u+1-i}\,.  \tag 5.13 
$$ 
Here and below we assume $\Re u\ll0$.  
 
On the other hand, by \tht{5.4} and \tht{5.9} we have 
$$ 
\gathered 
\pi_t(F^*(u))=1+\sum_{m=1}^\infty\frac{(t_1)_m(t_2)_m}{(-u)_mm!}\\ 
={}_2F_1(t_1,t_2;-u;1) 
=\frac{\Ga(-u)\Ga(-u-t_1-t_2)}{\Ga(-u-t_1)\Ga(-u-t_2)}\,,   
\endgathered 
\tag 5.14 
$$ 
where the last equality is Gauss' summation formula, cf. \tht{2.12}. When 
$t_1,t_2$ are as in \tht{5.11}, this coincides with \tht{5.13}. Thus, we have 
checked \tht{5.12} on the one--row functions.  
 
Next, we employ recurrence relations which make it possible to express 
the two--row functions through the one--row ones. It is convenient to 
extend the definition of two--row functions $P^*_{(p,q)}$ with 
$p>q\ge1$ to a larger set of indices by adopting the convention 
$P^*_{(p,q)}=-P^*_{(q,p)}$. Thus, $P^*_{(p,q)}$ makes sense for any $p,q=1,2,\dots$ (we assume that $P^*_{(p,q)}=0$ for $p=q$). Then we 
have two families of relations: 
$$ 
\gather 
P^*_{(p,1)}=P^*_{(p)}P^*_{(1)}-pP^*_{(p)}\,,  
\tag 5.15 
\\ 
P^*_{(p+1,q)}+P^*_{(p,q+1)}+(p+q)P^*_{(p,q)}= 
P^*_{(p)}P^*_{(q+1)}-P^*_{(p+1)}P^*_{(q)}-(p-q)P^*_{(p)}P^*_{(q)}\,.  
\tag 5.16 
\endgather 
$$ 
Here $p,q$ range over $\{1,2,\dots\}$. These relations were proved in 
\cite{I2}. Note that \tht{5.15} can be formally obtained from 
\tht{5.16} by substituting $q=0$.   
 
Using double 
induction on $p+q$ and $q$ we see that these relations indeed 
allow to express the two--row functions through the one--row 
functions.  
 
A direct computation shows that the relations \tht{5.15}, \tht{5.16} 
remain valid if we apply $\pi_t$ to each $P^*$ function involved. This means that \tht{5.12} holds on the two--row functions. 
 
Finally, to handle arbitrary $P^*$ functions we employ the following 
relation proved in \cite{I2}: 
$$ 
P^*_\mu 
=\Pf\left[P^*_{(\mu_i,\mu_j)}\right]_{1\le i,j\le\ell(\mu)+\varepsilon} 
\,. 
\tag 5.17 
$$ 
Here the symbol $\Pf$ means Pfaffian and $\varepsilon$ equals 1 for 
odd $\ell(\mu)$ and 0 for even $\ell(\mu)$, so that the order of the matrix 
is always even. Note that this relation has exactly the same form as 
in the case of the classical Schur $P$ functions, see \cite{M, III.8}.  
 
To conclude that \tht{5.12} holds on any $P^*_\mu$ we must verify that 
\tht{5.17} remains valid if we apply $\pi_t$ to each $P^*$ function. This 
readily follows from the well--known relation  
$$ 
\prod_{1\le i<j\le\ell(\mu)}\frac{\mu_i-\mu_j}{\mu_i+\mu_j} 
=\Pf\left[\frac{\mu_i-\mu_j} 
{\mu_i+\mu_j}\right]_{1\le i,j\le\ell(\mu)+\varepsilon}\,, \tag 5.18 
$$ 
see \cite{M, III.8}. \qed 
\enddemo 
 
\proclaim{Proposition 5.2} The harmonic function $\varphi_t$ afforded by 
Theorem 5.1 is strictly positive if and only if $t>0$.  
\endproclaim 

The proof is straightforward.
 
It should be noted that the values $t=k(1-k)/2$ used in the proof of 
Theorem 5.1 lie outside the region $t>0$.   
 
Quite similarly to the case of the Young graph, extreme points of 
$\Harm_1^+(\SS)$ correspond to {\it projective characters} or {\it 
projective finite factor representations} of the group $S(\infty)$. 
They can be parametrized by the points of the infinite--dimensional 
simplex $\Om(\SS)=\Om(\K)$ described by (4.7), see \cite{N}, \cite{I1}.   

The Poisson kernel $K(\mu,\omega)$ for the Schur graph is given by 
the image of the Schur $P$ function $P_\mu$ under the specialization 
of the algebra $\Gamma$ which sends odd powers sums to the following 
expressions (cf. (2.16), (3.9)) $$  
(x_1^k+x_2^k+\ldots)\mapsto\cases 1,& k=1,\\  
\sum_{i=1}^\infty\alpha_i^k,& k=2m+1\ge 3. 
\endcases 
\tag 5.19 
$$ 
\head \S6. Finite--dimensional specializations 
\endhead 
In the previous four sections we described four different examples of 
graphs which fit into the general scheme introduced in \S1. For each 
of these graphs we produced a nontrivial family of 
specializations of the corresponding algebras $A^*$ which defined, 
according to Proposition 1.3, a certain family of (nonnegative) 
harmonic functions on the graph.   
 
Every such family, in its turn, gives rise to a family of probability 
measures on the space $\Omega(\G)$ (Theorem 1.1), and this space in 
all our examples is infinite--dimensional, see (2.15), (4.7).  
 
For Young and Kingman graphs such measures have been thoroughly 
studied, see \cite{P.I--P.V}, \cite{BO1}, \cite{BO2}, \cite{Ki3}, \cite{Pi}, \cite{Ke4}. They lead to certain stochastic 
processes on the real line, for the Young graph the processes are 
closely related to those arising in Random Matrix Theory, while for 
the Kingman graph the theory is connected with Poisson processes.   
 
Our goal in this section is to construct `simpler' families of 
harmonic functions for Young, Kingman, and Schur graphs. The word 
`simpler' means that the corresponding measures on $\Omega(\G)$ will 
be supported by finite--dimensional subspaces. These measures will be 
explicitly computed.   
 
The corresponding Poisson integrals (which arise due to Theorem 1.1) 
will give (possibly new) integral formulas involving  products of Schur $S$ and $P$ functions, see 
\tht{6.10}, \tht{6.15}, \tht{6.23}, \tht{6.26} below.  
 
\subhead 6.1. Truncated Young branching \endsubhead Recall that for 
the Young graph $\Y$ the algebra $A^*$ is the algebra $\La^*$ of 
functions in infinitely many variables $x_1,x_2,\ldots$ symmetric in 
`shifted' variables $x_j'=x_j-j$.   
 
We shall consider the most natural specializations of this algebra 
obtained by fixing finitely many variables $x_1,\ldots,x_l$ and 
sending remaining variables $x_{l+1},x_{l+2},\ldots$ to zero.

In any such specialization all functions $s_\mu^*$ with the length (number of nonzero parts) of $\mu$ greater than $l$ vanish, see \cite{OO}. This implies that the harmonic function on $\Y$ afforded by Proposition 1.3 vanishes on all Young diagrams with more than $l$ rows. Thus, one can consider such a function as a harmonic function on the subgraph $\Y(l)$ of $\Y$ consisting of all Young diagrams with length $\le l$. This graph fits into the general formalism of $\S1$: the algebra $A$ is the algebra of symmetric polynomials in $l$ variables, and the algebra $A^*$ is the algebra of shifted symmetric polynomials in $l$ variables. The elements $\P_\mu$ and $\P_\mu^*$ are conventional and shifted Schur polynomials in $l$ variables, respectively. The graph $\Y(l)$ is called the {\it truncated Young graph}. Harmonic functions on such graphs were considered by Kerov, see \cite{Ke3}; he used them to derive certain Selberg--type integrals. Our arguments below are similar to those of Kerov's work.

Let us fix a Young diagram $\lambda$. Denote by $l$ the length of $\lambda$. We shall assume that $l\ge 2$. Denote by $\pi_\lambda$ the algebra homomorphism 
$\pi_\lambda:\Lambda^*\to\R$ defined by  
$$ 
\cases 
x_i\mapsto-\la_i-2(l-i)-1, &1\le i\le l,\\ 
x_i \mapsto 0 &i>l. 
\endcases 
\tag 6.1 
$$ 
 
Convenience of such choice of notation will be clear in a while. 
 
According to the general scheme of \S1 (Proposition 1.3), the 
corresponding harmonic function on $\Y$ has the form  
$$ 
\varphi_{\lambda}(\mu)= 
\frac{(-1)^{|\mu|}\pi_\la(s^*_\mu)}{(-\pi_\la(s^*_{(1)}))_{|\mu|}}= 
\frac{(-1)^{|\mu|}s^*_\mu(-\lambda-2\delta-1)}{(|\lambda|+l^2)_{|\mu|}}, 
\tag 6.2 
$$ 
where 
$s^*_\mu$ is the shifted Schur function and  
$$ 
\cases 
\delta_i=l-i, &1\le i\le l,\\ 
\delta_i=0, &i>l. 
\endcases 
$$  

\proclaim{Proposition 6.1} The function $\varphi_\lambda$ defined by 
\tht{6.2} is nonnegative.
\endproclaim 
\demo{Proof} Let the symbol $(a\dhr k)$ denote the {\it $k$th falling factorial 
power of $a$}:  
$$ 
(a\dhr k)=\cases a(a-1)\cdots (a-k+1),&k=1,2,\dots\,,\\ 1,&k=0. 
\endcases 
$$ 
Recall the definition of the shifted Schur polynomials in finitely many 
variables \cite{OO}  
$$ 
s_\mu^*(x_1,\ldots,x_l)= 
\frac{\det[(x_i+\delta_i\dhr \mu_j+\delta_j)]_{i,j=1}^l} 
{\det[(x_i+\delta_i\dhr \delta_j)]_{i,j=1}^l}\,, 
\tag 6.3 
$$ 
where $\delta_j$ are as above. 
 
Let us plug in our $x_i$ from \tht{6.1} to 
\tht{6.3}. We get  
$$ 
\gathered 
(-1)^{|\mu|}\pi_{\lambda}(s^*_\mu)= 
(-1)^{|\mu|}\, 
\frac{\det[(-\lambda_i-\delta_i-1\dhr \mu_j+\delta_j)]_{i,j=1}^l} 
{\det[(-\lambda_i-\delta_i-1\dhr \delta_j)]_{i,j=1}^l}\\= 
\frac{\det[\Gamma(\lambda_i+\delta_i+1+\mu_j+\delta_j)]_{i,j=1}^l} 
{\det[\Gamma(\lambda_i+\delta_i+1+\delta_j)]_{i,j=1}^l}. 
\endgathered 
$$

As the 
first $l$ members of sequences  
$\lambda+\delta+1$, $\mu+\delta$, and $\delta$ decrease,
the claim follows from the inequality
$$
\det[\Gamma(x_i+y_j)]_{i,j=1}^l>0, \qquad x_1>\dots>x_l> 0,\quad y_1>\dots>y_l> 0,  
$$
which is a special case of Problem VII.66 in \cite{PS}.  \qed
\enddemo

Next, we form the probability distributions $M_n=M_n^\la$ on $\Y_n$ for each $n=0,1,2,\dots$ according to (1.4): 
$$ 
M_n^\la(\nu)=\dim_\Y(\nu)\,\varphi_\la(\nu),\qquad |\nu|=n. 
\tag 6.4 
$$

Let us embed $\Y_n \hookrightarrow \Om(\Y)$ via 
$$ 
\nu\mapsto \left(\frac {\nu_1}{n},\frac{\nu_2}{n},\dots;\,0,0,\dots\right)
\in \Om(\Y).
\tag 6.5
$$ 
 
As was mentioned in \S1, the probability measure $P$ involved in the Poisson integral (1.2) is the weak limit of the images of the measures $M_n$ under appropriate embeddings $\Y_n \hookrightarrow \Om(\Y)$; as such  embeddings one can take (6.5).

Let us denote by $V(a)$ the Vandermonde determinant 
$$ 
V(a_1,\dots,a_s)=\prod_{1\le i<j\le s}(a_i-a_j). 
$$ 
 
\proclaim{Proposition 6.2}  
The images of $M_n^\la$ under the embeddings \tht{6.5} weakly converge to a probability measure $P$ on $\Om(\Y)$  supported by the 
finite--dimensional face  
$$ 
\gather 
\Delta_l= 
\{(\alpha,\beta)\in\Om(\Y)\,|\,\al_{l+1}= 
\al_{l+2}=\dots=\beta_1=\beta_2=\dots=0,\ \sum_{i=1}^l\al_i=1\}\\ \backsimeq 
\{(\al_1,\ldots,\al_l)\in \Bbb R^l_+|\,\al_1\ge 
\al_2\ge\dots\ge \al_l,\ \sum_{i=1}^l\al_i=1 \}. 
\endgather 
$$ 
The density of $P$ with respect to the Lebesgue measure on $\Delta_l$ equals  
$$ 
\frac{\Gamma(|\lambda|+l^2)\, 
s_\lambda(\al_1,\dots,\al_l)V^2(\al_1,\ldots,\al_l)} 
{\det[\Gamma(\lambda_i+\delta_i+\delta_j+1)]_{i,j=1}^l}. 
\tag 6.6 
$$  
\endproclaim 
 
\demo{Proof} We shall need the following lemma. 
\proclaim{Lemma 6.3} As $n\to\infty$, 
$$
\max_{\nu\in\Y_n}M_n^\la(\nu)=O(n^{1-l}).
$$
Furthermore, if $\nu_l\ge \varepsilon n$, where $\varepsilon>0$ is arbitrary, then, as $n\to\infty$,
$$ 
M_n^{\lambda}(\nu)= 
n^{1-l}\,\frac{\Gamma(|\lambda|+l^2)\, 
s_\lambda(\al_1,\dots,\al_l)V^2(\al_1,\ldots,\al_l)} 
{\det[\Gamma(\lambda_i+\delta_i+\delta_j+1)]_{i,j=1}^l}\,(1+o(1)) 
\tag 6.7 
$$
with
$$ 
\alpha_1=\frac{\nu_1}n,\dots, \alpha_l= 
\frac{\nu_l}n,
$$ 
and the estimate \tht{6.7} is uniform in $\nu$. 
\endproclaim 
 
Let us postpone the proof of this statement and proceed with the 
proof of Proposition 6.2 taking Lemma 6.3 for granted.

The Thoma simplex $\Om(\Y)$ defined in (2.15) is a compact topological space. Let $C(\Om(\Y))$ denote the algebra of continuous functions on $\Om(\Y)$. Take any $f\in C(\Om(\Y))$. Note that $f$ is bounded.  The value on $f$ of the image of $M_n^{\lambda}$ under the $n$th embedding (6.5) have the form  
$$ 
\sum_{\nu\in \Y_n}M_n^{\lambda}(\nu)\cdot 
f\left(\frac{\nu_1}n, \frac{\nu_2}n,\ldots;\,  
0,0,\ldots\right). 
\tag 6.8 
$$ 
 
Recall that $M_n^{\lambda}$ is supported by the Young diagrams $\nu$ with the length $\le l$. Let us first consider the part of 
the sum (6.8) involving diagrams $\nu$ with $\nu_l\ge \varepsilon n$, $\nu_{l+1}=0$. Using Lemma 6.3 and the boundedness of $f$, we get  
$$ 
\gathered 
\sum_{\Sb \nu\in \Y_n\\ 
\nu_l\ge\varepsilon n,\,\nu_{l+1}=0\endSb} M_n^{\lambda}(\nu)\cdot f\left(\frac{\nu_1}n, \frac{\nu_2}n,\ldots;\,  
0,0,\ldots\right)
= \frac{{\Gamma(|\lambda|+l^2)}\,(1+o(1))} 
{\det[\Gamma(\lambda_i+\delta_i+\delta_j+1)]_{i,j=1}^l}  
\\ \times\sum_{\Sb \nu\in \Y_n\\  
\nu_l\ge \varepsilon n,\, \nu_{l+1}=0\endSb}f(\al_1,\dots,\al_l,0,\dots;0,0,\dots)\cdot s_\lambda(\al)V^2(\al)\cdot n^{1-l}\,, 
\endgathered 
\tag 6.9 
$$ 
where $\al=(\al_1,\dots,\al_l)$ is as in Lemma 6.3. 
 
The sum in the right--hand side of (6.9) is a Riemannian sum for the 
integral  
$$ 
\int f(\al_1,\dots,\al_l,0,\dots;0,0,\dots)\cdot s_\lambda(\alpha)V^2(\alpha)\,d\alpha 
$$
over the part of $\Delta_l$ specified by the condition $\al_l\ge \varepsilon$.
 
Thus, it remains to prove that the part of the sum (6.8) involving 
diagrams $\nu$ with $\nu_l<\varepsilon n$ is $\varepsilon$--negligible.  
Since $M_n^\la(\,\cdot\,)=O(n^{1-l})$, this follows from the fact that the number of Young diagrams $\nu=(\nu_1,\dots,\nu_l,0,\dots)$ with $\nu_1+\dots+\nu_l=n$ and such that $\nu_l<\varepsilon n$ is bounded by $const\cdot \varepsilon\, n^{l-1}$. This completes the proof of Proposition 6.2 modulo Lemma 6.3. 
\enddemo 
 
\demo{Proof of Lemma 6.3} We shall employ formulas (6.2), (6.3), 
(6.4).
 Denote by $m$ the length of $\nu$. 
We apply a well--known dimension formula   
$$ 
\dim_{\Y}\nu=\frac{n!}{\prod_{i=1}^{m}(\nu_i+m-i)!} 
\prod_{1\le i<j\le m} (\nu_i-i-\nu_j+j) 
$$ 
which can be derived, e.g., from \cite{M, I.7, Ex.6}.
Then 
$$ 
\gathered 
M_n^{\la}(\nu)=\dim_{\Y}(\nu)\varphi_\la(\nu)= 
\frac{n!}{\prod_{i=1}^{m}(\nu_i+m-i)!} 
\prod_{1\le i<j\le m} (\nu_i-i-\nu_j+j)\\ \times 
\frac 1{(|\lambda|+l^2)_{n}} 
\cdot\frac{\det[\Gamma(\lambda_i+\delta_i+1+\nu_j+\delta_j)]_{i,j=1}^l} 
{\det[\Gamma(\lambda_i+\delta_i+1+\delta_j)]_{i,j=1}^l}. 
\endgathered 
$$ 
Here $\nu_j=0$ for $j>m$. 

Next, we have the following asymptotic relations as $n\to\infty$:
$$ 
\prod_{1\le i<j\le m} (\nu_i-i-\nu_j+j)= 
O(n^{\frac{m(m-1)}{2}}), 
$$ 
$$
\multline
\frac{\Gamma(\lambda_i+\delta_i+1+\nu_j+\delta_j)}{(\nu_j+m-j)!}= 
\frac{(\lambda_i+\delta_i+\nu_j+\delta_j)!}{(\nu_j+m-j)!}\\ \le (\lambda_i+\delta_i+\nu_j+\delta_j)^{\la_i+\delta_i+l-m} 
=O(n^{\la_i+\delta_i+l-m}), 
\endmultline
$$ 
$$ 
\frac{n!}{(|\lambda|+l^2)_{n}}= 
O(n^{-|\lambda|-l^2+1}), 
$$  
where estimates are uniform in $\nu\in \Y_n$.

Expanding the determinant $\det[\Gamma(\lambda_i+\delta_i+1+\nu_j+\delta_j)]_{i,j=1}^l$ and using the above estimates we get that each of $l!$ terms in the expansion, after multiplication by the remaining factors, is at most of order $n^{1-l}$. This proves the first claim of the lemma.

Let us proceed to the second claim. Assume that $\nu_l\ge \varepsilon n$. We have
$$ 
\prod_{1\le i<j\le l} (\nu_i-i-\nu_j+j)= 
n^{\frac{l(l-1)}{2}}\,V(\alpha)\cdot(1+o(1)), 
$$ 
$$ 
\frac{\Gamma(\lambda_i+\delta_i+1+\nu_j+\delta_j)}{(\nu_j+\delta_j)!}= 
\frac{\Gamma(\lambda_i+\delta_i+1+\nu_j+\delta_j)}{\Gamma(\nu_j+\delta_j+1)} 
=(\alpha_j n)^{\la_i+\delta_i}\cdot(1+o(1)), 
$$ 
$$ 
\frac{n!}{(|\lambda|+l^2)_{n}}= 
\Gamma(|\lambda|+l^2)\,\frac{\Gamma(n+1)}{\Gamma(|\lambda|+l^2+n)}= 
\Gamma(|\lambda|+l^2)\, n^{-|\lambda|-l^2+1}\cdot(1+o(1)). 
$$ 
All estimates are uniform in $\nu\in \Y_n$ provided that $\nu_l\ge \varepsilon n$. This assumption has been used in the second estimate above. This implies the desired estimate (6.7). 
\qed 
\enddemo 
 
As a consequence, we have the Poisson integral representation for our 
harmonic functions (see Theorem 1.1):  
\proclaim{Proposition 6.4}  
$$ 
\frac{(-1)^{|\mu|}s_{\mu}^*(-\la-2\delta-1)} 
{(|\la|+l^2)_{|\mu|}}=\frac{{\Gamma(|\lambda|+l^2)}} 
{\det[\Gamma(\lambda_i+\delta_i+\delta_j+1)]_{i,j=1}^l} 
\int_{\Delta_l}s_{\mu}(\alpha)\,s_{\lambda}(\alpha)\, V^2(\alpha)\,d\alpha. 
\tag 6.10 
$$ 
\endproclaim 
Note that the Poisson kernel $K(\mu,\om)$ on $\Delta_l$ coincides 
with the ordinary Schur function $s_{\mu}(\alpha)$, see the end of \S2.  
 
It should be noted that the integration in (6.10) can be carried out 
directly by making use of the formulas
$$
s_\mu(\al)=\frac{\det[\al_i^{\mu_j+\delta_j}]_{i,j=1}^l}{V(\al)},\quad
s_\la(\al)=\frac{\det[\al_i^{\la_j+\delta_j}]_{i,j=1}^l}{V(\al)},
$$ 
and the well--known  Dirichlet integrals
$$
\gathered
\int\limits_{\Sb\al_1\ge\dots\ge\al_l\ge 0\\ \al_1+\dots+\al_l=1\endSb}\al_1^{\kappa_1-1}\cdots\al_l^{\kappa_l-1}d\al =\frac 1{l!}\int\limits_{\Sb\al_1\ge 0,\dots,\al_l\ge 0\\ \al_1+\dots+\al_l=1\endSb}\al_1^{\kappa_1-1}\cdots\al_l^{\kappa_l-1}d\al\\ =\frac 1{l!}\, \frac{\Gamma(\kappa_1)\cdots\Gamma(\kappa_l)}{\Gamma(\kappa_1+\dots+\kappa_l)}.
\endgathered
$$
However, in more complicated cases, see below, a direct evaluation of integrals like (6.10) seems to be difficult. 

The identity (6.10) is similar to a result of Hua, see \cite{H, p. 104}, \cite{Ri, (3.3)}. 
 
\example{Remark 6.5} All claims of the present subsection remain 
true if we replace the integers $(\la_1\ge\la_2\ge\dots\ge\la_l>0)$ 
by any positive real numbers satisfying the same system of 
inequalities. Then the Schur function $s_\la(\al)$ should be understood 
just as the ratio  
${\det [\al_i^{\la_j+l-j}]}/{\det[\al_i^{l-j}]}.$ We restricted ourselves to integral $\la_i$'s in order to emphasize the symmetry $\la\leftrightarrow\mu$ in the integral (6.10). 
By analytic continuation, the formula (6.10) can be extended to arbitrary complex $\la_i$'s. When $\la$ has the form $((l-1)\theta+a,(l-2)\theta+a,\dots,a)$, where $\theta>0$ and $a>-1$, the measure (6.6) is related to the so--called Laguerre biorthogonal ensemble, see \cite{B2}.    
\endexample  
\subhead 6.2. $\Gamma$--shaped Young branching \endsubhead 
As in the previous subsection, we work with the Young graph $\Y$. 
This time we shall use another, so--called super realization of the algebra 
$A^*=\Lambda^*$ of shifted symmetric functions. Since a detailed exposition of the material of this subsection would be rather tedious, we shall only state the results and outline the ideas used in the proofs. 

Let $\wt \La$ be the algebra of supersymmetric functions in 
$x=(x_1,x_2,\dots)$ and $y=(y_1,y_2,\dots)$, see \cite{BR}, \cite{M, 
Ex.I.3.23-24, Ex.I.5.23}. It can be identified with the algebra $\La$ 
of symmetric functions; under this identification power sums 
$p_m\in\Lambda$ correspond to their super analogs  
$$ 
p_m(x;y)=\sum_ix_i^m+(-1)^{m-1}\sum_iy_i^m. 
$$   
Note that our notation slightly differs from that of Macdonald's 
book: his supersymmetric functions in $x$ and $y$ coincide with ours in $x$ and $-y$. 

Below we shall use Frobenius notation for Young diagrams, its description can be found in \cite{M, I.1}.
 
\proclaim{Theorem 6.6 \cite{KO}} There exists an algebra isomorphism 
$\rho:\La^*\to \wt \La$ such that for any $f\in\La^*$ and any Young 
diagram $\la=(\la_1,\la_2,\dots)$ with Frobenius coordinates 
$(p_1,\ldots p_d|\, q_1,\ldots,q_d)$ the following equality holds:  
$$ 
f(\la_1,\la_2,\dots)=\rho(f)\left({p_1+\frac 
12},\ldots,{p_d+\frac12};\, {q_1+\frac12},\ldots,{q_d+\frac12}\right). 
$$ 
\endproclaim 
 
Now we shall identify $\wt\La$ and $\La^*$ using the isomorphism $\rho$. We shall denote the elements $\rho(s^*_\mu)\in\wt\Lambda$ as $FS_\mu$ and call them {\it Frobenius--Schur functions}, see \cite{ORV}.

Let us consider 
specializations of the algebra $\La^*\simeq\La$ obtained by fixing  finitely many variables $x_1,\dots,x_d;\, y_1,\dots,y_d$ 
and sending remaining variables $x_i,y_i,\ i=d,d+1,\dots$, to zero. In any such specialization all functions $s^*_\mu$ with the depth (number of diagonal boxes) greater than $d$ vanish, this follows from results of \cite{ORV}. Hence, the corresponding harmonic functions are concentrated on a subgraph of the Young graph $\Y$ consisting of diagrams with depth $\le d$. These are exactly the Young diagrams which fit into the $\Gamma$--shaped figure with $d$ rows and $d$ columns. We denote the subgraph of such Young diagrams by $\Y(d,d)$ and call it the {\it $\Gamma$--shaped Young graph}. Like the truncated Young graphs considered in 6.1, the $\Gamma$--shaped Young graphs also fit into the general formalism of \S1. The algebras $A$ and $A^*$ are both identified with the algebra of supersymmetric polynomials in $d+d$ variables, the elements $\P_\mu$ are supersymmetric Schur polynomials, and the elements
$\P_\mu^*$ are supersymmetric Frobenius--Schur polynomials.  

Fix a Young diagram $\lambda$ with Frobenius coordinates 
$(p_1,\ldots p_d\,|\, q_1,\ldots,q_d)$. We shall denote by $\pi_\la$ 
the algebra homomorphism $\pi_\la:\La^*\to\R$ defined by 
$$ 
\cases 
x_i\mapsto -p_i-\frac 12,&1\le i\le d,\\ 
y_i\mapsto -q_i-\frac 12,&1\le i\le d,\\ 
x_i,\,y_i\mapsto 0,&i>d.
\endcases 
\tag 6.11 
$$   
 
According to \S1, the harmonic function on $\Y$ corresponding to 
$\pi_\la$ has the form  
$$ 
\varphi_\la(\mu)=\frac{(-1)^{|\mu|}\pi_\la(s^*_\mu)}{(-\pi_\la(s^*_{(1)}))_{|\mu|}}=\frac{(-1)^{|\mu|}FS_\mu(-p-\frac 12;-q-\frac 
12)} {(|\la|)_{|\mu|}}. 
\tag 6.12 
$$ 

Unfortunately, it is not clear how to prove directly that $\varphi_\la$ is nonnegative. However, we can go around this obstacle.

Set
$$ 
M_n^\la(\nu)=\dim_\Y(\nu)\,\varphi_\la(\nu),\qquad |\nu|=n.
$$
Consider the embeddings $\Y_n \hookrightarrow \Om(\Y)$ defined as follows. For a Young diagram $\nu\in\Y_n$ with Frobenius coordinates $(P_1,\dots,P_D|\,Q_1,\dots, Q_D)$ 
$$
\nu \mapsto \left(\frac{P_1+1/2}n,\dots,\frac{P_D+1/2}n,0,\dots;\,\frac{Q_1+1/2}n,\dots,\frac{Q_D+1/2}n,0,\dots\right)\in\Om(\Y).
\tag 6.13
$$

\proclaim{Proposition 6.7} The images of (possibly signed) measures $M_n^\la$ under the embeddings \tht{6.13} weakly converge, as $n\to\infty$, to a (positive) probability measure $P$ on $\Om(\Y)$. This measure is supported by the finite--dimensional face  
$$ 
\gather 
\Delta_{d,d}=\{(\alpha,\beta)\in\Om(\Y)\,|\,\al_{d+1}= 
\be_{d+1}=\al_{d+2}=\be_{d+2}=\dots=0,\ 
\sum_{i=1}^d(\al_i+\be_i)=1\}\\ \backsimeq 
\{(\al_1,\ldots,\al_d;\be_1,\dots,\be_d) 
\in \Bbb R^{2d}_+|\al_1\ge\dots\ge \al_d;\, 
\be_1\ge\dots\ge \be_d,\ \sum_{i=1}^d(\al_i+\be_i)=1 \}, 
\endgather 
$$ 
and its density with respect to the Lebesgue measure $d(\al;\be)$ on $\Delta_{d,d}$ equals  
$$
\frac{\Gamma(|\lambda|)} 
{\prod_{i=1}^d p_i!q_i!}\, 
\left[{\det}\left(\frac 1{p_i+q_j+1}\right)\right]^{-1}\, 
s_\lambda(\alpha;\beta) 
\,{\det}^2\left(\frac 1{\alpha_i+\beta_j}\right),
\tag 6.14 
$$ 
where $s_\la(\al;\be)$ is the supersymmetric Schur polynomial in $d+d$ variables. 
\endproclaim   

The proof of this proposition is quite similar to that of Proposition 6.2. An analog of Lemma 6.3 is proved using the Sergeev--Pragacz formula for $s_\la(\al;\be)$ (see \cite{PT}, \cite{M, I.3, Ex.24}) and its analog for the Frobenius--Schur polynomials (see \cite{ORV}). 

It turns out that Proposition 6.7 implies the existence of the Poisson integral representation (1.2) for the harmonic function (6.12). The proof of this claim is quite similar to the proof of Theorem B in \cite{KOO}. Since the Poisson kernel is always nonnegative, the existence of the Poisson integral representation implies that our harmonic function is nonnegative.  

Explicitly, the 
Poisson integral representation has 
the following form, cf. Proposition 6.4.  
 
\proclaim{Proposition 6.8} 
$$ 
\gathered 
\frac{(-1)^{|\mu|}FS_\mu(-p-\frac 12;-q-\frac 12)}{(|\la|)_{|\mu|}}= 
\frac{\Gamma(|\lambda|)}{\prod_{i=1}^d p_i!q_i!}\, 
\left[{\det}\left(\frac 1{p_i+q_j+1}\right)\right]^{-1}\\  
\times \int_{\Delta_{d,d}}s_{\mu}(\alpha;\be)\, 
s_{\lambda}(\alpha;\be)\,{\det}^2\left(\frac 
1{\alpha_i+\beta_j}\right) \, d(\alpha;\be). 
\endgathered 
\tag 6.15 
$$ 
\endproclaim 

\example{Remark 6.9} The formula (6.15) gives an expression for the integral 
$$
\int_{\Delta_{d,d}}s_{\mu}(\alpha;\be)\, 
s_{\lambda}(\alpha;\be)\,{\det}^2\left(\frac 
1{\alpha_i+\beta_j}\right) \, d(\alpha;\be).
\tag 6.16
$$
The integrand is symmetric in $\la$ and $\mu$. However, our assumptions on $\la$ and $\mu$ are different: the depth of $\la$ must be equal to $d$ while $\mu$ is arbitrary. Actually, if the depth of $\mu$ is $>d$ then both sides of (6.15) vanish. If the depth of $\mu$ is equal to $d$, the integration can be carried out directly in a rather simple way using the Berele--Regev formula
$$
s_\nu(\al_1,\dots,\al_d;\,\be_1,\dots,\be_d)={\left[\det\left(\frac 1{\al_i+\be_j}\right)_{i,j=1}^d\right]}^{-1}\,
\det[\al_i^{P_j}]_{i,j=1}^d\,\det[\be_i^{Q_j}]_{i,j=1}^d;
$$
here $(P|\,Q)$ are the Frobenius coordinates of $\nu$, see
\cite{BR}, \cite{M, I.3, Ex.23}. But if the depth of $\mu$ is strictly less than $d$, the Berele--Regev formula must be replaced by the more complicated Sergeev--Pragacz formula, and a direct integration seems to be more difficult.

\endexample

\example{Remark 6.10} All claims of this subsection remain true if 
we replace integral Frobenius coordinates $p_1>\dots>p_d\ge 0$, 
$q_1>\dots>q_d\ge 0$ of a fixed Young diagram $\la$ with 
any ordered sequences of real numbers $>-\frac 12$. Then the Schur 
function $s_\la(\al;\be)$ should be understood as 
$$
{\left[\det\left(\frac 1{\al_i+\be_j}\right)_{i,j=1}^d\right]}^{-1}\,
\det[\al_i^{p_j}]_{i,j=1}^d\,\det[\be_i^{q_j}]_{i,j=1}^d,
$$
cf. Remark 6.5. As in 6.1, we restricted ourselves to integral $p$'s and $q$'s in order to demonstrate the symmetry $\la\leftrightarrow\mu$ in (6.15). By analytic continuation, the formula (6.15) can be extrapolated to any pairwise distinct complex $p_i$'s and $q_i$'s.
\endexample

\subhead 6.3. Truncated Kingman branching \endsubhead
For the Kingman graph $\K$ (see \S4) the algebra $A^*$ coincides with
the algebra $\Lambda$ of symmetric functions. 
We consider specializations of $\Lambda$ obtained by fixing finitely
many indeterminates, say, $x_1,\dots,x_l$, and sending remaining indeterminates to zero. 

Under such a specialization all functions $\P_\mu^*=m_\mu^*$ with $\ell(\mu)>l$ vanish (this easily follows from the definition of the factorial monomial functions, see \S4). This means that the corresponding harmonic function lives on the subgraph $\K(l)$ of the Kingman graph $K$ consisting of all Young diagrams with the length $\le l$. We call this subgraph the {\it truncated Kingman graph}. It fits into the general formalism of \S1 with algebras $A$ and $A^*$ both equal to the algebra of symmetric polynomials in $l$ variables. The elements $\P_\mu$ and $\P_\mu^*$ are monomial symmetric polynomials and factorial monomial symmetric polynomials, respectively. 

Certain harmonic functions on truncated Kingman graphs and their applications to Selberg--type integrals were previously considered by Kerov \cite{Ke3}.

Let us fix a Young diagram 
$$
\la=(\la_1\ge\dots\ge\la_l>0)=(1^{r_1(\la)}2^{r_2(\la)}\dots).
$$
We define an algebra homomorphism $\pi_\la:\La\to \R$ as follows, cf.
(6.1), (6.11), 
$$
\cases
x_i\mapsto-\la_i-1, &1\le i\le l,\\
x_i\mapsto 0, &i>l.
\endcases
\tag 6.17
$$
The corresponding harmonic function on $\K$ has the form
$$
\varphi_{\lambda}(\mu)=\frac{(-1)^{|\mu|}\pi_\la(m^*_\mu)}{(-\pi_\la(m^*_{(1)}))_{|\mu|}}=
\frac{(-1)^{|\mu|}m^*_\mu(-\lambda_1-1,\dots,-\la_l-1)}{(|\lambda|+l)_{|\mu|}}.
\tag 6.18
$$  

\proclaim{Proposition 6.11} The harmonic function $\varphi_\la$ defined by \tht{6.18} is nonnegative. 
\endproclaim

\demo{Proof}
It suffices to note that
$m^*_\mu(-\lambda_1-1,\dots,-\la_l-1)$, by definition, is a sum of expressions of the
form ($\sigma$ is a permutation here) 
$$
\prod_{i=1}^{l(\mu)}(-\la_{\sigma(i)}-1)(-\la_{\sigma(i)}-2)
\cdots(-\la_{\sigma(i)}-\mu_i),
$$
each of which has sign $(-1)^{|\mu|}$.\qed
\enddemo

According to \S1, we have probability distributions $M_n^\la$ on $\K_n$
for all $n=1,2,\ldots$ given by 
$$
M_n^{\la}(\nu)=\dim_\K(\nu)\,\varphi_{\la}(\nu), \quad |\nu|=n.
\tag 6.19
$$
Embeddings $\K_n\hookrightarrow\Om(\K)$ are defined as follows
($\Om(\K)$ was defined in (4.7)): 
$$
\nu\in\K_n\mapsto\left(\frac{\nu_1}n,\frac{\nu_2}n,\dots\right)\in\Om(\K).
\tag 6.20
$$
The images of the probability distributions $M_n^{\la}$ via these
embeddings, according to the general theory, weakly converge, as $n\to\infty$, to a certain
probability measure $P$ on $\Om(\K)$. 
\proclaim{Proposition 6.12} The probability measure $P$ on
$\Om(\K)$ is supported by the finite--dimensional face 
$$
\gather
\Delta_l=\{\alpha\in\Om(\K)\,|\,\al_{l+1}=\al_{l+2}=\dots=0,\ 
\sum_{i=1}^l\al_i=1\}\\ \backsimeq
\{(\al_1,\ldots,\al_l)\in \Bbb R^l_+|\,\al_1\ge
\al_2\ge\dots\ge \al_l,\ \sum_{i=1}^l\al_i=1 \}.
\endgather
$$
Its density with respect to the Lebesgue measure $d\al$ on $\Delta_l$ equals 
$$
\frac{\Gamma(|\la|+l)
\cdot r_1(\la)!r_2(\la)!\cdots}{\prod_{i=1}^l\la_i!}\,
m_\la(\al_1,\dots,\al_l).
\tag 6.21
$$ 
\endproclaim

The proof is very similar to that of Proposition 6.2, so we shall
just state the analog of Lemma 6.3 in this case. 
\proclaim{Lemma 6.13} As $n\to\infty$, 
$$
\max_{\nu\in\K_n} M_n^\la(\nu)=O(n^{1-l}).
$$
Furthermore, if $\nu_{l+1}=0$ and $\nu_l\ge \varepsilon n$, where $\varepsilon>0$ is arbitrary, then, as $n\to\infty$, 
$$
M_n^{\la}(\nu)=n^{1-l}\, \frac{\Gamma(|\la|+l)
\cdot r_1(\la)!r_2(\la)!\cdots}
{\prod_{i=1}^l\la_i!}\,m_\la(\al_1,\dots,\al_l)\,(1+o(1)),
\tag 6.22
$$
with $
\al_1={\nu_1}/n,\dots,\al_l={\nu_l}/n.
$
The estimate \tht{6.22} is uniform in $\nu$.
\endproclaim
\demo{Proof} Using the formula (4.6) for $\dim_\K(\nu)$ we get
$$
M_n^{\la}(\nu)=\dim_\K(\nu)\,\varphi_\la(\nu)=
\frac{n!}{(|\la|+l)_{n}}\,\frac {\sum_{\sigma} 
\prod_{i=1}^{\ell(\nu)}(\la_{\sigma(i)}+1)_{\nu_i}}{\nu_1!\cdots\nu_{\ell(\nu)}!},
$$
where the summation is taken over all permutations $\sigma\in S_l$ which produce different products $\prod_{i=1}^{\ell(\nu)}(\la_{\sigma(i)}+1)_{\nu_i}$.

Next, we have asymptotic relations
$$
\frac {n!}{(|\la|+l)_n}=O(n^{1-l-|\la|}),
$$
$$
\frac{(\la_j+1)_{\nu_i}}{\nu_i!}\le (\la_j+\nu_i)^{\la_j}=O(n^{\la_j}).
$$
They imply that each term of the sum above, after multiplication by remaining factors, is at most of order
$n^{1-l}$. This proves the first part of the lemma.

For the second part of the lemma, assume $\nu_l\ge \varepsilon n$. Then we have   
$$
\frac {n!}{(|\la|+l)_n}=\Gamma(|\la|+l)\, n^{1-l-|\la|}(1+o(1)),
$$
$$
\frac{(\la_j+1)_{\nu_i}}{\nu_i!}=\frac {(\al_in)^{\la_j}}{\la_j!}\,(1+o(1))
$$
as $n\to\infty$, all estimates are uniform in $\nu\in \K_n$ provided that $\nu_l\ge \varepsilon n$. This yields the estimate (6.22).\qed 
\enddemo

As a corollary, we get the Poisson integral representation, cf.
Propositions 6.4 and 6.8. 

\proclaim{Proposition 6.14}
$$
\gathered
\frac{(-1)^{|\mu|}m^*_\mu(-\lambda_1-1,\dots,-\la_l-1)}{(|\lambda|+l)_{|\mu|}}\\=
\frac{\Gamma(|\la|+l)\cdot r_1(\la)!r_2(\la)!\cdots}{\prod_{i=1}^l
\la_i!}\int_{\Delta_l}m_\mu(\al)m_\la(\al)d\al.
\endgathered
\tag 6.23
$$
\endproclaim

This claim, similarly to Proposition 6.4, can be proved directly by making use of Dirichlet integrals. 
\example{Remark 6.15} All claims above remain valid for any positive ordered sequence $\la=(\la_1\ge\la_2\ge\dots\ge\la_l\ge 0)$ with the
obvious modification of the definition of monomial symmetric function $m_\la(\al)$, cf. Remarks 6.5, 6.10. By analytic continuation, the formula (6.23) can be extended to arbitrary complex $\la_i$'s.  
\endexample

\subhead 6.4. Truncated Schur branching \endsubhead
In this subsection we shall deal with the Schur graph, see \S5. For this graph the algebras $A$ and $A^*$ coincide with a subalgebra $\Gamma$ of the algebra of symmetric functions; $\Gamma$ is generated by the odd power sums $\sum_ix_i^{2k+1}$, $k=0,1,\dots$ . Again, we consider specializations of $\Gamma$ obtained by fixing variables $x_1,\dots,x_l$ and sending remaining variables to zero. As in 6.2, we shall state the results and sketch the ideas of the proofs.  

In such a specialization, the elements $\P_\mu^*=P_\mu^*$ vanish if $\ell(\mu)>l$. This means that the corresponding harmonic function can be viewed as a harmonic function on the {\it truncated Schur graph} $\SS(l)$ --- the subgraph of the Schur graph $\SS$ consisting of diagrams with length $\le l$. The truncated Schur graphs also fit into the formalism of \S1 with algebras $A$ and $A^*$ coinciding with the subalgebra of the algebra of symmetric polynomials in $l$ variables generated by odd power sums. The elements $\P_\mu$ and $\P_\mu^*$ are the Schur $P$ polynomials and the factorial Schur $P$ polynomials, respectively.

Let us fix a strict partition $\lambda$ and denote its length by $l$.
We define a multiplicative linear functional 
$\pi_\la:\Gamma\to \R$ as follows, cf. (6.1), (6.11), (6.17):
$$
\cases
x_i\mapsto-\la_i-1, &1\le i\le l,\\
x_i\mapsto 0, &i>l.
\endcases
$$
The corresponding harmonic function on $\SS$ has the form
$$
\varphi_{\lambda}(\mu)=\frac{(-1)^{|\mu|}\pi_\la(P^*_\mu)}{(-\pi_\la(P^*_{(1)}))_{|\mu|}}=\frac{(-1)^{|\mu|}P^*_\mu(-\lambda_1-1,\dots,-\la_l-1)}
{(|\lambda|+l)_{|\mu|}}.
\tag 6.24
$$

As in 6.2, it is not evident that this function is nonnegative.

For all strict partitions $\nu\in\SS_n$, $n=1,2,\ldots$, we define, as usual,
$$
M_n^{\la}(\nu)=\dim_\SS\nu\cdot\varphi_\la(\nu).
$$
Embeddings $\SS_n \hookrightarrow\Om(\SS)=\Om(\K)$ are defined
exactly as in the case of the Kingman graph, see (6.20), the only
difference is that now all partitions are strict.  

\proclaim{Proposition 6.16} The sequence of images of (possibly signed) measures $M_n^{\la}$ under the embeddings defined above weakly
converges, as $n\to\infty$, to a (positive) probability measure $P$ on $\Om(\SS)$. This
measure is supported by the finite--dimensional face 
$$
\gather
\Delta_l=\{\alpha\in\Om(\SS)\,|\,\al_{l+1}=\al_{l+2}=\dots=0,\ 
\sum_{i=1}^l\al_i=1\}\\ \backsimeq
\{(\al_1,\ldots,\al_l)\in \Bbb R^l_+|\,\al_1\ge
\al_2\ge\dots\ge \al_l,\ \sum_{i=1}^l\al_i=1 \},
\endgather
$$ 
and its density with respect to the Lebesgue measure $d\al$ on $\Delta_l$ equals 
$$
\frac{\Gamma(|\lambda|+l)}{\prod_{i=1}^l\lambda_i!}
\left[\pf\left(\frac{\lambda_i-\lambda_j}{\lambda_i+\lambda_j+2}\right)
\right]^{-1}P_\lambda(\alpha_1,\dots,\al_l)\,
{\pf}^2\left(\frac{\alpha_i-\alpha_j}{\alpha_i+\alpha_j}\right).
\tag 6.25
$$
\endproclaim 

This result is parallel to Propositions 6.2, 6.7, 6.12. Its proof is based on an appropriate analog of the approximation Lemmas 6.3, 6.13. The proof of such a lemma in this case follows from explicit formulas for Schur $P$--functions and factorial Schur functions \cite{M, III.8}, \cite{I1}.  

Similarly to Proposition 6.7, Proposition 6.16 implies the existence of the Poisson integral representation for $\varphi_\la$. Thanks to the positivity of (6.25), this implies that $\varphi_\la$ is nonnegative. 
 
The explicit Poisson integral representation (1.2) in this case takes the following form.

\proclaim{Proposition 6.17}
$$
\gathered
\frac{(-1)^{|\mu|}P^*_\mu(-\lambda_1-1,\dots,-\la_l-1)}{(|\lambda|+l)_{|\mu|}}=\frac{\Gamma(|\lambda|+l)}
{\prod_{i=1}^l\lambda_i!}
\left[\pf\left(\frac{\lambda_i-\lambda_j}{\lambda_i+\lambda_j+2}\right)
\right]^{-1}\\ \times
\int\limits_{\Delta_l}P_\mu(\al)P_\la(\al)\pf^2\left(\frac{\al_i-\al_j}
{\al_i+\al_j}\right)\,d\al.
\endgathered
\tag 6.26
$$
\endproclaim

\example{Remark 6.18} If $\ell(\mu)=l\;(=\ell(\la))$, then the integration in (6.26) can be carried out directly using Dirichlet integrals and the formula \cite{M, III.8, Ex.12}
$$
P_\nu(\al_1,\dots,\al_l)=\left[\pf\left(\frac{\al_i-\al_j}{\al_i+\al_j}\right)\right]^{-1}
\det[\al_i^{\nu_j}]_{i,j=1}^l\,,
\tag 6.27
$$
which holds for all $\nu$ of length $l$. If $\ell(\mu)<l$, the integration seems to be more complicated, cf. Remark 6.9. 
\endexample

\example{Remark 6.19} All claims of the present subsection remain true for any nonnegative strictly ordered sequence
$\la=(\la_1>\la_2>\dots>\la_l\ge 0)$. Then the Schur $P$--function
$P_\la(\al)$ should be understood as 
$$
\left[\pf\left(\frac{\al_i-\al_j}{\al_i+\al_j}\right)\right]^{-1}
\det[\al_i^{\la_j}]\,,
$$ 
cf. Remarks
6.5, 6.10, 6.15.  
By analytic continuation, the formula (6.26) can be extended to arbitrary complex mutually distinct $\la_i$'s.
\endexample
 
\head \S7. Appendix \endhead 
 
\demo{Proof of Theorem 1.1} Existence of the integral representation 
\tht{1.2} follows from Choquet's theorem, see, e.g., \cite{Ph}. To 
prove its uniqueness one can apply another theorem, due to Choquet 
and Meyer, \cite{DM}. Then we have to verify that the cone 
$\Harm^+(\G)$ is a lattice, i.e., for any 
$\varphi,\psi\in\Harm^+(\G)$, there exist their lowest upper bound 
$\varphi\vee\psi$ and greatest lower bound $\varphi\wedge\psi$. Let 
us prove that 
$$ 
\gather 
(\varphi\vee\psi)=\lim_{n\to\infty}\sum_{\la\in\G_n}  
\dimG(\mu,\la)\max(\varphi(\la), \psi(\la)),   \tag7.1   \\ 
(\varphi\wedge\psi)=\lim_{n\to\infty}\sum_{\la\in\G_n}  
\dimG(\mu,\la)\min(\varphi(\la), \psi(\la)).   \tag7.2   \\ 
\endgather 
$$ 
 
Indeed, take \tht{7.1} and \tht{7.2} as the definition of the 
functions $\varphi\vee\psi$ and $\varphi\wedge\psi$. For any fixed 
$\mu$, the sum in the right--hand side of \tht{7.1}  
increases as $n\to\infty$ and remains bounded from above by 
$\varphi(\mu)+\psi(\mu)$. Similarly, the sum in the right--hand side of \tht{7.2} 
decreases and remains bounded from below by 0. Hence, the limits 
exist.  
 
Next, remark that for any fixed vertex $\nu$, the function 
$\mu\mapsto\dimG(\mu,\nu)$ satisfies the harmonicity relation up to 
level $|\nu|-1$. This implies that $\varphi\vee\psi$ and 
$\varphi\wedge\psi$ are harmonic functions. Clearly, they are 
nonnegative and are upper and lower bounds, respectively.  
 
Finally, it is readily verified that they are the lowest upper bound and 
the greatest lower bound, respectively. \qed 
\enddemo

\bigskip 
\bigskip 
 
{\smc A.~Borodin}: Department of Mathematics, The University of 
Pennsylvania, Philadelphia, PA 19104-6395, U.S.A.   
 
E-mail address: 
{\tt borodine\@math.upenn.edu} 
 
{\smc G.~Olshanski}: Dobrushin Mathematics Laboratory, Institute for 
Problems of Information Transmission, Bolshoy Karetny 19, 101447 
Moscow GSP-4, RUSSIA.   
 
E-mail address: {\tt olsh\@iitp.ru, 
olsh\@glasnet.ru} 
 
\enddocument 
 

\Refs  
\widestnumber\key{KOO} 
 
\ref\key BR 
\by A.~Berele and A.~Regev 
\paper Hook Young diagrams with applications to combinatorics and to representations of Lie superalgebras 
\jour Adv. Math. 
\vol 64 
\yr 1987 
\pages 118--175 
\endref 
 
\ref\key B1 
\by A.~M.~Borodin 
\paper Multiplicative central measures on the Schur graph 
\inbook Representation theory, dynamical systems, combinatorial and 
algorithmical methods II (A.~M.~Vershik, ed.) 
\bookinfo Zapiski Nauchnykh Seminarov POMI {\bf 240} 
\publ Nauka 
\publaddr St.~Petersburg 
\yr 1997 
\pages 44--52 (Russian) 
\transl English transl. to appear in J. Math. Sci. 
\endref  
 
\ref\key B2 
\bysame 
\paper Biorthogonal ensembles 
\jour Nucl. Phys. B 
\vol 536 
\yr 1998 
\pages 704--732 (preprint version available via  
{\tt http://xxx.lanl.gov/abs/math/9804027})  
\endref 
 
\ref\key BO1 
\by A.~Borodin and G.~Olshanski 
\paper Point processes and the infinite symmetric group 
\jour Math. Research Lett. 
\vol 5
\yr 1998
\pages 799--816 (preprint version available via 
{\tt http://xxx.lanl.gov/abs/ math/9810015}) 
\endref
 

\ref\key BO2
\bysame
\paper Distributions on partitions, point processes and the hypergeometric kernel
\paperinfo Pre\-print, 1999, available via 
{\tt http://xxx.lanl.gov/abs/math/9904010} 
\endref

\ref\key DM
\by C.~Dellacherie, P.-A.~Meyer
\book Probabilities and potential
\publ North--Holland Elsevier
\yr 1978, 1982
\endref


\ref\key Er 
\by A.~Erdelyi (ed.)  
\book Higher transcendental functions, {\rm Vol. 1} 
\publ Mc Graw--Hill 
\yr 1953 
\endref 
 
\ref\key Ew 
\by W.~J.~Ewens 
\paper Population Genetics Theory -- the Past and the Future 
\inbook  Mathematical and statistical developments of 
evolutionary theory.  Proc. NATO ASI Symp. 
\ed S.~ Lessard 
\publ Kluwer 
\publaddr Dordrecht 
\yr 1990 
\pages 117--228 
\endref 
 
\ref \key H
\by L.-K.~Hua
\book Harmonic analysis of functions of several complex variables in the classical domains
\bookinfo Translations of Math. Monographs, Vol. 6
\publ Amer. Math. Soc.
\publaddr Providence, R.I.
\yr 1963
\endref

\ref\key I1 
\by V.~N.~Ivanov 
\paper Dimension of skew shifted Young diagrams and 
projective representations of the infinite symmetric group 
\inbook Representation theory, dynamical systems, combinatorial and 
algorithmical methods II (A.~M.~Vershik, ed.) 
\bookinfo Zapiski Nauchnykh Seminarov POMI {\bf 240} 
\publ Nauka 
\publaddr St.~Petersburg 
\yr 1997 
\pages 115--135 (Russian) 
\transl\nofrills English transl. to appear in J. Math. Sci. 
\endref 
 
\ref\key I2 
\bysame 
\paper Paper in preparation 
\endref 
 
\ref\key JK 
\by G.~James, A.~Kerber 
\book The representation theory of the symmetric group 
\bookinfo Encyclopedia of mathematics and its applications {\bf 16} 
\publ Addison--Wesley 
\yr 1981 
\endref 
 
\ref\key Ke1 
\by S.~V.~Kerov 
\paper Combinatorial examples in the theory of AF-algebras 
\inbook Differential geometry, Lie groups and mechanics X 
\bookinfo Zapiski Nauchnykh Seminarov LOMI, Vol. 172 
\yr 1989 
\pages 55-67  
\lang Russian 
\transl English translation in J. Soviet Math. {\bf 59} (1992), 
No.~5, pp. 1063--1071 
\endref 
 
\ref\key Ke2 
\bysame 
\paper Generalized Hall--Littlewood symmetric functions and orthogonal 
polynomials 
\inbook Representation Theory and Dynamical Systems 
\ed A.~M.~Vershik 
\bookinfo Advances in Soviet Math. {\bf 9} 
\publ Amer. Math. Soc. 
\publaddr Providence, R.I. 
\yr 1992 
\pages 67--94 
\endref 
 
\ref\key Ke3 
\bysame 
\paper The boundary of Young lattice and random Young  
tableaux 
\inbook Formal power series and algebraic combinatorics (New Brunswick,  
NJ, 1994) 
\bookinfo DIMACS Ser. Discrete Math. Theoret. Comput. Sci. 
\vol 24 
\yr 1996 
\pages 133--158 
\publ Amer. Math. Soc. 
\publaddr Providence, RI 
\endref 
  
\ref\key Ke4 
\bysame 
\paper Subordinators and permutation actions with quasi-invariant 
measure 
\inbook Representation theory, dynamical systems, combinatorial and 
algorithmical methods I (A.~M.~Vershik, ed.) 
\bookinfo Zapiski Nauchnykh Seminarov POMI {\bf 223} 
\publ Nauka 
\publaddr St.~Petersburg 
\yr 1995 
\pages 181--218 
\lang Russian 
\transl English transl. in J. Math. Sciences {\bf 87} (1997), No. 6 
\endref 
 
\ref\key Ke5 
\bysame 
\paper Anisotropic Young diagrams and Jack symmetric functions 
\jour Funct. Anal. Appl. 
\pages to appear (preprint version available via  
{\tt http://xxx.lanl.gov/abs/math/9712267}) 
\endref 
 
\ref \key KO  
\by S.~Kerov and G.~Olshanski 
\paper Polynomial functions on the set of Young diagrams 
\jour Comptes Rendus Acad.\ Sci.\ Paris S\'er. I 
\vol 319 
\yr 1994 
\pages 121--126 
\endref 
 
\ref\key KOO 
\by S.~Kerov, A.~Okounkov, G.~Olshanski 
\paper The boundary of Young graph with Jack edge multiplicities 
\jour Intern. Math. Res. Notices  
\yr 1998 
\pages  No.4, 173--199 
\endref 
 
\ref \key KOV  
\by S.~Kerov, G.~Olshanski, A.~Vershik  
\paper Harmonic analysis on the infinite symmetric group. A deformation  of the regular representation  
\jour Comptes Rend. Acad. Sci. Paris, S\'er. I  
\vol 316  
\yr 1993  
\pages 773--778  
\endref  
 
\ref\key KV1 
\by S.~V.~Kerov and A.~M~Vershik  
\paper Characters, factor representations and $K$-functor of the 
infinite symmetric group.  
\inbook Operator algebras and group representations, Vol. II (Neptun, 
1980) 
\bookinfo Monographs Stud. Math. {\bf 18}  
\publ Pitman, Boston--London 
\yr 1984 
\pages 23--32 
\endref 
 
\ref\key KV2  
\by S.~Kerov, A.~Vershik  
\paper The Grothendieck group of the infinite symmetric group and  
symmetric functions with the elements of the $K_0$-functor theory  
of AF-algebras  
\inbook Representation of Lie groups and related topics  
\bookinfo Adv. Stud. Contemp. Math. {\bf 7}  
\eds A.~M.~Vershik and D.~P.~Zhelobenko  
\publ Gordon and Breach  
\yr 1990  
\pages 36--114  
\endref  
 
\ref\key Ki1 
\by J.~F.~C.~Kingman 
\paper Random partitions in population genetics 
\jour Proc. Roy. Soc. London A. 
\vol 361 
\yr 1978 
\pages 1--20 
\endref 
 
\ref\key Ki2 
\bysame 
\paper The representation of  partition structures  
\jour J. London Math. Soc. 
\vol 18 
\yr 1978 
\pages 374--380 
\endref 
 
\ref\key Ki3 
\bysame 
\book Poisson processes 
\publ Oxford University Press 
\yr 1993 
\endref 
 
\ref\key KS 
\by F.~Knop and S.~Sahi 
\paper Difference equations and symmetric polynomials defined by 
their zeros 
\jour Internat. Math. Res. Notices 
\yr 1996 
\issue 10 
\pages 473--486 
\endref 
 
\ref\key M  
\by I.~G.~Macdonald  
\book Symmetric functions and Hall polynomials  
\bookinfo 2nd edition  
\publ Oxford University Press  
\yr 1995  
\endref 
 
\ref\key N
\by M.~Nazarov
\paper Projective representations of the infinite symmetric group
\inbook
Representation theory and dynamical systems (A.~M.~Vershik, ed.), Adv. Soviet Math., 9, Amer. Math. Soc., Providence,
RI
\pages 115--130
\yr 1992
\endref


\ref 
\key Ok1 
\bysame 
\paper 
Quantum immanants and higher Capelli identities 
\jour Transformation Groups 
\vol 1 \issue 1-2 \yr 1996 \pages 99--126 
\endref 
 
\ref 
\key Ok2 
\bysame 
\paper 
(Shifted) Macdonald polynomials: $q$-Integral 
representation and combinatorial formula 
\jour Comp. Math.
\vol 12
\yr 1998
\pages 147--182
\endref 
 
\ref 
\key Ok3 
\bysame  
\paper Thoma's theorem and representations of infinite bisymmetric
group 
\jour Funct. Anal. Appl. 
\vol 28
\yr 1994
\pages no. 2, 101--107
\endref

\ref\key Ok4
\bysame
\paper Binomial formula for Macdonald polynomials
\jour Math. Res. Lett.
\vol 4
\yr 1997
\pages 533--553
\endref

\ref 
\key OO 
\by A.~Okounkov and G.~Olshanski 
\paper Shifted Schur functions 
\jour Algebra i Analiz 
\vol 9 
\issue 2 
\yr 1997 
\pages 73--146 
\lang Russian 
\transl\nofrills English translation: St.~Petersburg Math. J.  
{\bf 9} (1998), no.~2, 239--300. 
\endref 
 
\ref 
\key OO2 
\bysame 
\paper Shifted Jack polynomials, binomial formula, 
and applications 
\jour Math.\ Res.\ Lett.\  
\vol 4 \yr 1997 \pages 69--78 
\paperinfo q-alg/9608020 
\endref 
 
\ref\key OV
\by A.~Okounkov and A.~Vershik
\paper A new approach to representation theory of symmetric groups
\jour Selecta Math. (N.S.)
\vol 2
\issue 4
\yr 1996
\pages 581--605
\endref


\ref\key Ol
\by G.~I.~Olshanskii
\paper Unitary representations of $(G,K)$-pairs 
connected with the infinite symmetric group $S(\infty)$ 
\jour Algebra i Analiz
\vol 1
\issue 4
\yr 1989
\pages 178--209
\lang Russian
\transl English translation in Leningrad Math.\ J. 
\vol 1 \issue 4 \pages 983-1014\yr 1990
\endref


\ref\key ORV 
\by G.~Olshanski, A.~Regev, and A.~Vershik 
\paper Frobenius--Schur functions
\paperinfo Preprint, 1999, to appear in ``xxx.lanl.gov'' archive  
\endref 
 
\ref\key  P.I 
\by G.~Olshanski 
\paper Point processes and the infinite symmetric group. Part I: The 
general formalism and the density function 
\paperinfo Preprint, 1998, available via  
{\tt http://xxx.lanl.gov/abs/ math/9804086} 
\endref 
 
\ref\key P.II 
\by A.~Borodin 
\paper Point processes and the infinite symmetric group. Part II: 
Higher correlation functions 
\paperinfo Preprint, 
 1998, available via 
 {\tt http://xxx.lanl.gov/abs/math/9804087} 
\endref 
 
\ref\key P.III 
\by A.~Borodin and G.~Olshanski 
\paper Point processes and the infinite symmetric group. Part III: 
Fermion point processes 
\paperinfo Preprint, 1998, available via  
{\tt http://xxx.lanl.gov/abs/math/ 9804088} 
\endref 
 
\ref\key P.IV 
\by A.~Borodin 
\paper Point processes and the infinite symmetric group. Part IV: 
Matrix Whittaker kernel 
\paperinfo Preprint, 
 1998, available via   
{\tt http://xxx.lanl.gov/abs/math/ 9810013} 
\endref 
 
\ref\key P.V 
\by G.~Olshanski 
\paper Point processes and the infinite symmetric group. Part V: 
Analysis of the matrix Whittaker kernel  
\paperinfo Preprint, 1998, available via  
{\tt http://xxx.lanl.gov/abs/math/ 9810014} 
\endref 
 
\ref\key Ph
\by R.~R.~Phelps
\book Lectures on Choquet's theorem 
\publ Van Nostrand mathematical studies, Vol. 7
\publaddr Princeton, NJ
\yr 1966
\endref

\ref\key Pi 
\by J.~Pitman 
\paper The two--parameter generalization of Ewens' random partition 
structure  
\paperinfo Univ. Calif. Berkeley, Dept. Stat. Technical Report {\bf 
345}  
\yr 1992 
\pages 1--23 
\endref 
 
\ref\key PY 
\by J.~Pitman and M.~Yor 
\paper The two--parameter Poisson--Dirichlet distribution derived 
from a stable subordinator 
\jour Ann. Prob. 
\vol 25 
\yr 1997 
\pages 855--900 
\endref 

\ref \key PS
\by G.~P\'olya and G.~Szeg\"o
\book Problems and theorems in analysis, Vol. II
\publ Springer--Verlag
\publaddr Berlin, etc.
\yr 1976
\endref
 
 
\ref\key Pr 
\by P.~Pragacz 
\paper Algebro-geometric applications of Schur $S$- and $Q$-polynomials 
\inbook Topics in Invariant Theory 
\bookinfo Seminaire d'Alg\`ebre Paul Dubriel et Marie-Paule 
Malliavin, Lecture Notes in Math. 
\vol 1478 
\publ Springer--Verlag 
\publaddr New--York/Berlin 
\yr 1991 
\pages 130--191 
\endref 
 
\ref\key PT 
\by P.~Pragacz and A.~Thorup 
\paper On a Jacobi--Trudi identity for supersymmetric polynomials 
\jour Adv. Math. 
\yr 1992f 
\pages 8--17 
\endref 
 
\ref\key Ri
\by D.~St.~P.~Richards
\paper Analogs and extensions of Selberg's integral
\inbook q-Series and partitions (D.~Stanton, ed.)
\publ Springer--Verlag
\publaddr New York, etc.
\yr 1989
\pages 109--137
\endref


\ref\key Ro 
\by N.~A.~Rozhkovskaya 
\paper Multiplicative distributions on Young graph 
\inbook Representation theory, dynamical systems, combinatorial and 
algorithmical methods II (A.~M.~Vershik, ed.) 
\bookinfo Zapiski Nauchnykh Seminarov POMI {\bf 240} 
\publ Nauka 
\publaddr St.~Petersburg 
\yr 1997 
\pages 246--257 (Russian) 
\transl\nofrills English transl. to appear in J. Math. Sci. 
\endref  
 
\ref 
\key S 
\by S.~Sahi 
\paper The spectrum of certain invariant differential operators 
associated to a Hermitian symmetric space 
\inbook Lie Theory and Geometry: In Honor of Bertram Kostant 
\eds J.-L.~Brylinski, R. Brylinski, V.~Guillemin, V. Kac 
\bookinfo Progress in Mathematics {\bf 123} 
\publ Birkh\"auser 
\publaddr Boston, Basel 
\yr 1994 
\pages 569--576 
\endref 
 
\ref\key T 
\by E.~Thoma 
\paper Die unzerlegbaren, positive--definiten Klassenfunktionen 
der abz\"ahlbar unendlichen, symmetrischen Gruppe 
\jour Math.~Zeitschr. 
\vol 85 
\yr 1964 
\pages 40-61 
\endref 
 
\ref\key VK 
\by A.~M.~Vershik, S.~V.~Kerov 
\paper Asymptotic theory of characters of the symmetric group 
\jour Funct. Anal. Appl.  
\vol 15 
\yr 1981 
\pages 246--255 
\endref 
 
\ref\key W
\by A.~J.~Wassermann
\paper Automorphic actions of compact groups on operator algebras
\paperinfo Thesis, University of Pennsylvania
\yr 1981
\endref

\endRefs
\bye